\documentclass{amsart}

\usepackage{amssymb,amscd, amsmath, amsfonts}
\usepackage[arrow,matrix]{xy}
\usepackage[mathscr]{euscript}
\CompileMatrices
\DeclareMathAlphabet{\mathssbx}{OT1}{cmss}{bx}{n} 
\begin{document}
 
\setcounter{tocdepth}{1}

\newcommand\datver[1]{\def\datverp%
 {\par\boxed{\boxed{\text{Version: #1; Run: \today}}}}}

\font\rsfself=rsfs11  
\font\rsfsacht=rsfs8  
\font\rsfsfuenf=rsfs5 

\def\SS{{\mathchoice{{\mbox{{\rsfself S}\hspace*{3pt}}}}%
                   {{\mbox{{\rsfself S}\hspace*{3pt}}}}%
                   {{\mbox{{\rsfsacht S}\hspace*{3pt}}}}%
                   {{\mbox{{\rsfsfuenf S}\hspace*{1pt}}}}}}%

\renewcommand{\theenumi}{\alph{enumi}}
\renewcommand{\labelenumi}{\rm {({\theenumi})}}
\renewcommand{\labelenumii}{(\roman{enumii})}
%
%
\newcommand{\wt}[1]{\widetilde{#1}}
\newcommand{\Prod}{\prod}
\newcommand{\Cal}{\mathcal}
\newcommand\Ran{\operatorname{Im}}
\newcommand\abs{\operatorname{abs}}
\newcommand\rel{\operatorname{rel}}
\newcommand\inv{\operatorname{inv}}
\newcommand\topo{\operatorname{top}}
\newcommand\opp{\operatorname{op}}
\newcommand\mfk{\mathfrak}
\newcommand\coker{\operatorname{coker}}
\newcommand\hotimes{\hat \otimes}
\newcommand\ind{\operatorname{ind}}
\newcommand\End{\operatorname{End}}
\newcommand\per{\operatorname{per}}
\newcommand\pa{\partial}
\newcommand\sign{\operatorname{sign}}
\newcommand\supp{\operatorname{supp}}
\newcommand\cy{\mathcal{C}^\infty}
\newcommand\CI{\mathcal{C}^\infty}
\newcommand\CO{\mathcal{C}_0}
\newcommand\lra{\longrightarrow}
\newcommand\vlra{-\!\!\!-\!\!\!-\!\!\!\!\longrightarrow}
\newcommand\bS{{}^b\kern-1pt S}
\newcommand\bT{{}^b\kern-1pt T}
\newcommand\Hom{\operatorname{Hom}}

\newcommand\alg[1]{\mathfrak{A}(#1)} 
\newcommand\qalg[1]{\mathfrak{B}(#1)} 
\newcommand\ralg[1]{\mathfrak{A}_r(#1)} 
\newcommand\rqalg[2]{\mathfrak{B}_r(#1)}
\newcommand\ideal[1]{C^*(#1)} 
\newcommand\rideal[1]{C^*_r(#1)} 
\newcommand\qideal[2]{\mathfrak{R}_{#1}(#2)} 
\newcommand\In{\operatorname{In}}

\newcommand\TR{\operatorname{T}}
\newcommand\ha{\frac12} 
\newcommand\cal{\mathcal}
\newcommand\END{\operatorname{END}}
\newcommand\ENDG{\END_{\GR}(E)}
\newcommand\CC{\mathbb C}
\newcommand\NN{\mathbb N}
\newcommand\RR{\mathbb R}
\newcommand\R{\mathbb R}
\newcommand\ZZ{\mathbb Z}
\newcommand\ci{${\mathcal C}^{\infty}$}
\newcommand\CIc{{\mathcal C}^{\infty}_{\text{c}}}
\newcommand\hden{{\Omega^{\lambda}_d}}
\newcommand\VD{{\mathcal D}}
\newcommand\mhden{{\Omega^{-1/2}_d}}
\newcommand\ehden{r^*(E)\otimes {\Omega^{\lambda}_d}}


\newcommand{\Cat}{\mathcal C}
\newcommand{\Gr}[1]{{\mathcal G}^{(#1)}}
\newcommand{\GR}{\mathcal G}
\newcommand{\LGR}{\mathcal L}
\newcommand{\BB}{\mathbb{B}}
\newcommand{\GG}{\mathcal G}
\newcommand{\OA}{\mathcal O}
\newcommand{\tPS}[1]{\Psi^{#1,0}(\GR)}
\newcommand{\AL}{{\mathcal A}(\GR)}
\newcommand{\FAM}{P=(P_x,x \in \Gr0)}
\newcommand\symb[2]{{\mathcal S}^{#1}(#2)}
\newcommand{\loc}{\operatorname{loc}}
\newcommand{\cl}{\operatorname{cl}}
\newcommand{\A}{s}
\newcommand{\prop}{\operatorname{prop}}
\newcommand{\comp}{\operatorname{comp}}
\newcommand{\adb}{\operatorname{adb}}
\newcommand{\dist}{\operatorname{dist}}

\newcommand{\alp}{r }
\newcommand{\bet}{d }
\newcommand{\gm}{\Gamma }
\newcommand{\lon}{\longrightarrow }
\newcommand{\be}{\begin{eqnarray*}}
\newcommand{\ee}{\end{eqnarray*}}
\newcommand{\GGR}{{\GR}}
\newcommand{\cald}{{\cal D}}
\newcommand{\calx}{{\cal X}}
\def\cinfo{{\mathcal C}^{\infty,0}}
\def\ccinfo{{\mathcal C}_c^{\infty,0}}
\def\ccinf0{{\mathcal C}_c^{\infty,0}}
\newcommand{\II}{\ccinfo(S^*(\GR), \End(E) \otimes {\mathcal P}_m)}

\newcommand{\IIY}{\CIc(S^*(\GR\vert_Y), 
\End(E\vert_Y) \otimes {\mathcal P}_m)}
\newcommand{\mI}{\mathfrak I}

\newcommand{\cfg}{continuous family groupoid}

\def\nin{\noindent}
\def\eg{e.g.\ }
\def\pt#1#2{{\partial #1\over \partial #2}}

\newcommand{\frakg}{{\mathfrak g}}

\let\Tilde=\widetilde
\let\Bar=\overline
\let\Vec=\overrightarrow
\let\ceV=\overleftarrow
\def\vlra{\hbox{$\,-\!\!\!-\!\!\!-\!\!\!-\!\!\!-\!\!\!
-\!\!\!-\!\!\!-\!\!\!-\!\!\!-\!\!\!\longrightarrow\,$}}

\def\vleq{\hbox{$\,=\!\!\!=\!\!\!=\!\!\!=\!\!\!=\!\!\!
=\!\!\!=\!\!\!=\!\!\!=\!\!\!=\!\!\!=\!\!\!=\!\!\!=\!\!\!=\,$}}

\def\lrah{\hbox{$\,-\!\!\!-\!\!\!
-\!\!\!-\!\!\!-\!\!\!-\!\!\!-\!\!\!\longrightarrow\,$}}

\def\surj{-\!\!\!-\!\!\!-\!\!\!\gg}

\def\inj{>\!\!\!-\!\!\!-\!\!\!-\!\!\!>}
  

\newcommand{\ad}[1]{{\mathcal A}^{\fd}_{#1}}
\newcommand{\as}{\wt{a}}
\newcommand{\aun}{\mathcal{A}_{\infty}}
\newcommand{\aund}{\ad{\infty}}
\newcommand{\bo}{{}^{b}\Omega}
\newcommand{\bu}{{\mathcal B}_{U}}
\newcommand{\bs}{\wt{b}}
\newcommand{\ca}{{\mathcal A}}
\newcommand{\cae}{\ca_{1}}
\newcommand{\can}{\ca_{0}}
\newcommand{\caun}{\ca_{\infty}}
\newcommand{\cav}{\ca_{\Omega}}
\newcommand{\cavs}{\wt{\ca}_{\Omega}}
\newcommand{\cb}{{\mathcal B}}
\newcommand{\cH}{{\mathcal{H}}}
\newcommand{\cj}{{\mathcal J}}
\newcommand{\cjt}[1]{\cj_{#1}(\ft)}
\newcommand{\cjtk}{\cj_{k}(\ft)}
\newcommand{\cjtu}{\cjt{\infty}} 
\newcommand{\ck}{{\mathcal K}}
\newcommand{\ckx}{{\mathcal K}_{x}}
\newcommand{\cnb}{\mathcal{C}_{b}}
\newcommand{\cun}{\mathcal{C}^{\infty}}
\newcommand{\cunb}{\mathcal{C}^{\infty}_{b}}
\newcommand{\cunc}{\cun_{c}}
\newcommand{\cuncd}{\dot{\mathcal{C}}^{\infty}_{c}}
\newcommand{\cv}{\mathcal{V}}
\newcommand{\cz}{\CC}
\newcommand{\dbar}{d\hspace{-2.8pt}\rule[5.5pt]{3pt}{0.24pt}}
\newcommand{\dd}{\Delta_{\fd}}
\newcommand{\de}[1]{{\mathcal D}(#1)}
\newcommand{\deh}{{\mathcal D}^{\frac{1}{2}}}
\newcommand{\dez}{\delta_{12}}
\newcommand{\dirint}{\int^{\oplus}}
\newcommand{\djx}{D_{j;x}}
\newcommand{\djxs}{\wt{D}_{j;x}}
\newcommand{\dmy}{D_{m;y}}
\newcommand{\djys}{\wt{D}_{j;y}}
\newcommand{\dmjx}{\delta_{M_{j,x}}}
\newcommand{\dt}{\delta_{T}}
\newcommand{\dxs}{\frac{dx'}{x'}}
\newcommand{\dxss}{dx''}
\newcommand{\elm}{E_{\ell,m}}
\newcommand{\fd}{\mathfrak{D}}
\newcommand{\fm}{\mathfrak{M}}
\newcommand{\fs}{\mathfrak{S}}
\newcommand{\ft}{\mathfrak{T}}
\newcommand{\gpt}{\mathcal{G}_{\Psi}(T)}
\newcommand{\gx}{{\GR_{x}}}
\newcommand{\hol}[1]{{\cal O}(#1)}
\newcommand{\hft}[1]{{\mathcal H}^{#1}(\ft)}
\newcommand{\huft}{\hft{\infty}}
\newcommand{\idh}{{\rm id}_{H}}
\newcommand{\lch}{\mathcal{L}(\cH)}
\newcommand{\lck}{\mathcal{L}(\ck)}
\newcommand{\lckx}{\mathcal{L}(\ckx)}
\newcommand{\ldch}{\mathcal{L}_{D}(\cH)}
\newcommand{\ldgch}{\mathcal{L}^{\GR}_{D}(\cH)}
\newcommand{\lh}{\mathcal{L}(H)}
\newcommand{\lhx}{\mathcal{L}(H_{x})}
\newcommand{\llzrny}{\mathcal{L}(\lzrny)}
\newcommand{\lukdh}{{\mathcal L}^{\infty}_{K,d}(H)}
\newcommand{\lukh}{{\mathcal L}^{\infty}_{K}(H)}
\newcommand{\lzrny}{L^{2}(\rny)}
\newcommand{\lzgx}{L^{2}(\gx;\rsed)}
\newcommand{\mjx}{M_{j;x}}
\newcommand{\mmy}{M_{m;y}}
\newcommand{\norm}[1]{\|#1\|}
\newcommand{\nz}{\NN}
\newcommand{\olt}{\omega_{T}^{\ell}}
\newcommand{\olta}[1]{\omega_{#1}^{\ell}}
\newcommand{\olrt}{\omega_{T_{1},T_{2}}^{\ell,r}}
\newcommand{\ort}{\omega_{T}^{r}}
\newcommand{\orta}[1]{\omega_{#1}^{r}}
\newcommand{\ot}{\omega_{T}}
\newcommand{\otimesh}{\wh{\otimes}}
\newcommand{\pege}{\Psi^{1,0}(\GR;\deh)}
\newcommand{\pft}[1]{\Psi^{#1}(\ft)}
\newcommand{\pmge}{\Psi^{m,0}(\GR;\deh)}
\newcommand{\pmunge}{\Psi^{-\infty,0}(\GR)}
\newcommand{\pnge}{\Psi^{0,0}(\GR)}
\newcommand{\pmungee}{\Psi^{-\infty,0}(\GR)}
\newcommand{\pngee}{\Psi^{0,0}(\GR)}
\newcommand{\puft}{\pft{\infty}}
\newcommand{\re}{R^{(1)}}
\newcommand{\rklx}{\RR^{k,\ell}_{x}}
\newcommand{\rkx}{\RR^{k}_{x}}
\newcommand{\rne}{\RR^{n}_{\eta}}
\newcommand{\rny}{\RR^{n}_{y}}
\newcommand{\rpq}{\overline{\rz}_{+}}
\newcommand{\rsed}{r^{*}\deh}
\newcommand{\rz}{\RR}
\newcommand{\rj}{R^{(j)}}
\newcommand{\sk}{\Sigma_{K}}
\newcommand{\spk}[1]{\langle #1 \rangle}
\newcommand{\sva}{\sigma_{\Omega}(a;\varphi_{0},\psi_{0})}
\newcommand{\te}{T^{(1)}}
\newcommand{\tz}{T^{(2)}}
\newcommand{\tj}{T^{(j)}}
\newcommand{\wh}[1]{\widehat{#1}}


%
%
\newcommand\Mand{\text{ and }}
\newcommand\Mandset{\text{ and set}}
\newcommand\Mas{\text{ as }}
\newcommand\Mat{\text{ at }}
\newcommand\Mfor{\text{ for }}
\newcommand\Mif{\text{ if }}
\newcommand\Min{\text{ in }}
\newcommand\Mon{\text{ on }}
\newcommand\Motwi{\text{ otherwise }}
\newcommand\Mwith{\text{ with }}
\newcommand\Mwhere{\text{ where }} 
\newcommand\ie{{\em i.e., }} 
%
%
\newtheorem{theorem}{Theorem}
\newtheorem{proposition}{Proposition}
\newtheorem{corollary}{Corollary}
\newtheorem{lemma}{Lemma}
\newtheorem{definition}{Definition}
\newtheorem{notation}{Notations}
\newtheorem{remark}{Remark}
\theoremstyle{remark}

%
%
\def\n#1#2{\| #1 \|_{{#2}}}
\def\g0{\GR^{(0)}}
\def\S{\mathscr{S}}
\def\cred#1{C^*_{\mathrm{r}}(#1)}
\def\rep#1{{\R_+^*}\!^{#1}}  
\def\r+#1{\R_+^{#1}}
\def\cinfo{{\mathcal C}^{\infty,0}}
\def\ccinfo{{\mathcal C}_c^{\infty,0}}
\def\Cc{{\mathcal C}_c}
\let \mx \mbox
\let \hx \hbox
\let \vx \vbox

\def\C{\mathbb{C}}
\def\E{\mathbb{E}}
\def\A{\mathbb{A}}
\def\B{\mathbb{B}}
\def\D{\mathbb{D}}
\def\F{\mathbb{F}}
\def\H{\mathbb{H}}
\def\N{\mathbb{N}}
\def\P{\mathbb{P}}
\def\Q{\mathbb{Q}}
\def\R{\mathbb{R}}
\def\Z{\mathbb{Z}}
\def\G{\mathbb{G}}
\def\K{\mathbb{K}}
\def\J{\mathbb{J}}
\def\T{\mathbb{T}}

\def\diagr#1{\def\normalbaselines{\baselineskip=0pt
\lineskip=10pt\lineskiplimit=1pt} \matrix{#1}}
\def\hfl#1#2{\smash{\mathop{\hx to 12mm{\rightarrowfill}}
\limits^{\scriptstyle#1}_{\scriptstyle#2}}}
\def\vfl#1#2{\llap{$\scriptstyle #1$}\left\downarrow
\vx to 6mm{}\right.\rlap{$\scriptstyle #2$}}
\def\antihfl#1#2{\smash{\mathop{\hx to 12mm{\leftarrowfill}}
\limits^{\scriptstyle#1}_{\scriptstyle#2}}}
\def\antivfl#1#2{\llap{$\scriptstyle #1$}\left\uparrow
\vx to 6mm{}\right.\rlap{$\scriptstyle #2$}}

\def\build#1_#2^#3{\mathrel{\mathop{\kern 0pt#1}\limits_{#2}^{#3}}}

\def\limind{\mathop{\oalign{lim\cr \hidewidth$\longrightarrow$\hidewidth\cr}}}
\def\limproj{\mathop{\oalign{lim\cr \hidewidth$\longleftarrow$\hidewidth\cr}}}

\def\norme#1{\left\| #1 \right\|}
\def\module#1{\left| #1 \right|}
\def\va#1{\left| #1 \right|}
\def\scal#1{\left\langle #1 \right\rangle}
\def\inv#1{{#1}^{-1}}
\def\invf{f^{-1}}

\def\g0{G^{(0)}}
\def\psd{pseudodiff\'erentiel}
\def\dg{\partial G}
\def\dm{\partial M}
\def\dx{\partial X}

\def\intm{\overset{\:\bullet}{M}}
\def\intx{\overset{\:\bullet}{X}}
\def\intf#1{\overset{\:\bullet}{#1}}
\def\cc#1{C^*(#1)}
\def\kc#1{K_*(C^*(#1))}
\def\rep#1{{\R_+^*}^{#1}}  
\def\r+#1{\R_+^{#1}}

\def\tg{groupo\"\i de tangent}
\def\gr{groupoid}
\def\d{{\rm d}}
\def\e{{\varepsilon}}
\def\cstar{$C^*$-al\-g\`ebre}
\def\rbar{\overline\R}
\def\l{\lambda}
\def\kth{$K$-th\'e\-o\-rie}
\def\ie{{\it i.e.}}
\def\cinf{$C^\infty$}
\def\M{M}
\def\G{{\bf G}}
\def\E{\mathcal{E}}
\def\In{\hbox{In}}
\def\pd{pseudodifferential}

\def\k{\mathssbx{k}}
\def\K{\mathcal{K}}
\def\I{\mathcal{I}}
\def\L{\mathcal{L}}

\def\lb{\hbox{lb}}
\def\rb{\hbox{rb}}

\def\bsp{$b$-stretched product}
\def\bcalc{$b$-calculus}

\def\sm{submanifold}
\def\mc{manifold with corners}
\def\mcs{manifolds with corners}

\def\D{\hbox{D}}
\def\S{\mathcal S}

\def\A{\mathcal A}
\def\F{\mathcal F}
\def\HF{\mathcal{HF}}
\def\sym{\mathfrak S}
\def\GG{\mathcal G}

\def\sh{\mathrm{sh}}
\def\ch{\mathrm{ch}}

\author[R. Lauter]{Robert Lauter} \address{Universit\"at 
       Mainz. Fachbereich 17-Mathematik, D-55099 Mainz, Germany} 
         \email{lauter@mathematik.uni-mainz.de}

\author[B. Monthubert]{Bertrand Monthubert}
       \address{Laboratoire Emile Picard, Universit\'e Paul Sabatier (UFR MIG),
         118 route de Narbonne, F-31062 Toulouse CEDEX 4 }
       \email{monthube@picard.ups-tlse.fr}

\author[V. Nistor]{Victor Nistor} \address{Pennsylvania State
       University, Math. Dept., University Park, PA 16802}
       \email{nistor@math.psu.edu} 

\thanks{Lauter was partly supported by a scholarship of the  
        German Academic Exchange Service (DAAD) within 
        the {\em Hochschulsonderprogramm III von Bund und 
        L\"andern}.
        Monthubert was partially supported by a ACI Jeunes Chercheurs and 
        by the European Network 
         ``Geometric Analysis'' (HPRN-CT-1999-00118).
       Nistor was partially supported by an  NSF Young
       Investigator Award DMS-9457859 and a Sloan Research
       Fellowship, and NSF Grant DMS-9971951. Manuscripts available from {\bf
       http:{\scriptsize//}www.math.psu.edu{\scriptsize/}nistor{\scriptsize/}}.}

\begin{abstract} 
We construct algebras of pseudodifferential operators on a continuous
family groupoid $\GR$ that are closed under holomorphic functional
calculus, contain the algebra of all pseudodifferential operators of
order $0$ on $\GR$ as a dense subalgebra, and reflect the smooth
structure of the groupoid $\GR$, when $\GR$ is smooth. As an
application, we get a better understanding on the structure of
inverses of elliptic pseudodifferential operators on classes of
non-compact manifolds.  For the construction of these algebras closed under
holomorphic functional calculus, we develop three methods: one
using 
semi-ideals, one using commutators, and one based on
Schwartz spaces on the groupoid. 
\end{abstract}

\title[Spectral invariance of algebras of pseudodifferential operators]
      {Spectral invariance for certain algebras of 
       pseudodifferential operators} 
      
\maketitle 
\tableofcontents

\def\c{\cite}
\def\fr{\frac}
\def\ub{\underbar}
\def\O{\Cal O}
\def\F{\Cal F}
\def\differ{\text{differentiable} }
\def\tPSeudo{\text{pseudodifferential} }
\def\supp{\text{supp} }
\def\inn{{\mathcal R}} 
\def\prop{\text{prop}} 
\def\frag{\frak{G}}
\def\simd{\tilde{d}}
\def\simf{\tilde{\F}}
\def\simo{\tilde{\O}}
\def\simr{\tilde{r}}
\def\simp{\tilde{p}}
\def\simmu{\tilde{\mu}} 
\def\O12{\Omega^\frac{1}{2}}

\section*{Introduction}

Let $M$ be a compact manifold and $P$ an elliptic pseudodifferential
operator of order $m \ge 0$ on $M$. Assume that $P$ is invertible as
an unbounded operator on $L^2(M)$ (the space of square integrable
$1/2$-densities on $M$). A classical and very useful result states
that then $P^{-1}$ is also a pseudodifferential operator.
For non-compact manifolds the situation is more complicated,
essentially because we also want to control the behavior at infinity
of the inverse.

For example, when $M$ has cylindrical ends, a convenient class of
pseudodifferential operators is that of $b$-pseudodifferential
operators introduced by Melrose \cite{MelroseACTA,Melrose42} (see also
\cite{rs86}). Then it is known that the inverse $P^{-1}$ of an
elliptic $b$-pseudodifferential operator (defined in $L^2$-sense) is
not necessarily also a $b$-pseudodifferential operator
in the so-called {\em small $b$-calculus}. We say that the
$b$-calculus is not {\em spectrally invariant}. There exist however
different classes of pseudodifferential operators associated to a
manifold with cylindrical ends that are spectrally invariant. In this
paper we discuss this property for the $c_n$-calculi, which are
spectrally invariant for $n \ge 2$.\footnote{When $n = 1$, the
$c_n$-calculus is nothing but the $b$-calculus and, when $n = 2$, it
is usually called the ``cusp-calculus,''
see for instance \cite{meni96c} (which is based on earlier work of
Melrose). Here, $n$ should not be confused with the dimension of the
manifold but determines the degree of degeneracy in direction to the
boundary, more precisely, the $c_n$-calculus is modeled on the
differential operators with degeneracies of the form $x^n\pa_x$ at the
boundary where $x$ stands for the direction normal to the boundary.}

Closely related to spectral invariance is the question whether a given
algebra of pseudodifferential operators, say of order $0$, is closed
under holomorphic functional calculus.  Let us explain the relevance
of this property.  Let $H$ be a Hilbert space and $a\in\lh$ {\em
relatively invertible} in $\lh$, i.e.\ there exists $\as\in\lh$ with
$a\as a=a$ and $\as a\as=\as$. By a characterization of Atkinson
\cite{atk}, we know that this is the case if, and only if, the range
$R(a)$ of $a$ is closed. By a classical result of Rickart
\cite{rick46} (see also \cite{cor67} and \cite[Bemerkung 5.7]{gra84}),
$0$ is an isolated point of the spectrum $\sigma(a^{*}a)$ of $a^{*}a$,
and the orthogonal projection $p$ onto the kernel $N(a^{*}a)=N(a)$ is
given by the integral
\begin{equation}
        \label{proj} p:=\frac{1}{2\pi
        i}\int_{\gamma}(z\idh-a^{*}a)^{-1}dz,
\end{equation}
where $\gamma$ is a small circle around the origin in $\cz$ that
does not intersect $\sigma(a^{*}a)$. In that case, 
the operator 
\begin{equation}
\label{relinv}
\as=(p+a^{*}a)^{-1}a^{*}
\end{equation}
is a relative inverse of $a$. Because of $(\as a)^{*}=\as a$ and
$(a\as)^{*}=a\as$ the operator $\as$ is also called an {\em orthogonal
generalized inverse} or a {\em Moore-Penrose inverse of $a$}
\cite{nashed}; it is easily seen to be uniquely determined.  If
$a:H\rightarrow H$ is additionally a Fredholm operator, then we have
$p=\idh-\as a$, and $q:=\idh-a\as$ is the orthogonal projection onto
the orthogonal complement of $R(a)$, i.e.\ $\as$ is a Fredholm inverse
of $a$. For more about generalized inverses we refer to \cite{nashed}.

In applications, we often know much more about the operator $a$ than
having closed range or being Fredholm, so it is natural to ask, which
properties of $a$ are inherited by the Moore-Penrose inverse $\as$.
It is immediate from \eqref{proj} and \eqref{relinv} that $\as$
belongs to the intersection of all subalgebras $A\subseteq\lh$ that
are symmetric with respect to the $*$-operation from $\lh$ and closed
under holomorphic functional calculus in $\lh$; in particular, any
property of $a$ that can be covered by a symmetric subalgebra $A$ of
$\lh$ that is closed under holomorphic functional calculus is true for
the Moore-Penrose inverse as well.  Thus, it is interesting to find
algebras that are closed under holomorphic functional calculus.
Without loss of generality, we can always assume that an algebra that
is closed under holomorphic functional calculus is also symmetric.  

Of particular importance for pseudodifferential and microlocal
analysis are symmetric, continuously embedded {\em Fr\'echet}
subalgebras of $C^{*}$-algebras that are closed under holomorphic
functional calculus ($\Psi^{*}$-algebras \cite{gra84}).  Indeed, in
contrast to the rather rigid $C^{*}$-topology, the Fr\'echet-topology
allows a flexible treatment of $\cun$-phenomena within a functional
analytic setting \cite{grapot92,gra98}.  On the other hand, stability
under holomorphic functional calculus and symmetry still establishes a
strong relation between the structure of a $\Psi^{*}$-algebra and that
of its $C^{*}$-closures leading to sometimes unexpected insights into
the internal structure of a $\Psi^{*}$-algebra; for instance, the set
of relatively invertible elements in $\Psi^{*}$-algebras has been shown
to be a locally rational Fr\'echet manifold \cite{gra84}.  Starting
from the seminal work \cite{gra84} of Gramsch, the world of
$\Psi^{*}$-algebras has been explored by many authors, and we refer
the reader to \cite{gra84,gra90,gra98,guw,jung95,aa,lor95} and the
references given there for more details.  Besides, the $K$-theory
groups, $K_{*}(A)$, of a symmetric subalgebra of a $C^{*}$-algebra
$\cb$ that is closed under holomorphic functional calculus coincide
with those of its closure in $\cb$.  Our basic example for an algebra
that is closed under holomorphic functional calculus is the algebra of
classical pseudodifferential operators of order $0$ on a closed
manifold.

In \cite{lmn} the authors considered a pseudodifferential calculus on
continuous family groupoids; this calculus generalizes the
pseudodifferential calculus on $\cinfo$-foliations used by Connes to
prove the index theorem for foliated spaces \cite{connesF}.  

See \cite{Taylor, TaylorB} for an introduction to the theory of
pseudodifferential operators. In a slightly different context, Nistor,
Weinstein, Xu \cite{nwx} and Monthubert, Pierrot \cite{mopi97} have
studied a pseudodifferential calculus on differentiable groupoids. As
demonstrated by the examples in \cite{mont97,nwx} and the survey
\cite{sur}, the groupoid approach yields a pseudodifferential calculus
for many interesting situations in analysis and geometry, especially
on open manifolds and manifolds with singularities, in a unified
way. Up to some support condition, this pseudodifferential calculus
recovers the existing calculi, and in many cases goes beyond the
results that are known in the literature.  Let us only mention that
the class of algebras that can be defined using groupoids include the
ordinary pseudodifferential calculus, a $G$-equivariant
pseudodifferential calculus on bundles of Lie groups, the $b$-calculus
of Melrose and many of its cousins on manifolds with corners, the edge
calculus on manifolds with fibered boundaries, the calculus of
adiabatic pseudodifferential operators, and many others.  On the other
hand, we know for many of the different pseudodifferential calculi
mentioned above that the algebra of operators of order $0$ is not
closed under holomorphic functional calculus. In fact, due to the
support condition in the case of a general continuous family groupoid
(which is a quite convenient condition that will insure that the
composition is defined), the algebra of operators of order $0$ is
almost never closed under holomorphic functional calculus.

In the present paper, we develop a general strategy to embed the
algebra $\pngee$ of pseudodifferential operators of order $0$ on a
continuous family groupoid $\GR$ into larger algebras $\ca$ that are closed
under holomorphic functional calculus and still share some of the
interesting properties with the algebra $\pngee$.  (We shall denote by
$\Psi^{m,0}(\GR)$ the space of order $m$ pseudodifferential operators
on a continuous family groupoid $\GR$.) In fact, it is one of the
results of this paper that it usually suffices to embed the algebra
$\pmunge$ of operators of order $-\infty$ in an algebra $\cj$ that is
closed under holomorphic functional calculus.  Up to some technical
conditions, $\ca:=\pnge+\cj$ is then an algebra that is closed under
holomorphic functional calculus.  For the construction of the algebra
$\cj$, we suggest 
three alternatives. The first one relates properties of an algebra to
those of a two-sided ideal and its corresponding quotient, the second
one is based on commutator methods from operator theory, whereas the
third one, depends more on the geometry of the groupoid and requires
the existence of a length function $\phi$ with polynomial growth on
the groupoid. The role of $\phi$ is to define a Schwartz space $\SS(\GR)$
on $\GR$ that replaces $\cunc(\GR)=\pmunge$. In this way, we control
the behavior at infinity of the kernels of our pseudodifferential
operators on a typical leaf $d^{-1}(x)$ of the groupoid.

The paper is organized as follows: In {\sc Section} \ref{Sec.1} we
recall the notion of algebras closed under functional calculus,
and consider
and answer 
the question, whether an algebra is closed under functional
calculus provided an ideal and the corresponding quotients are. In
{\sc Section} \ref{Sec.2} we develop the operator theoretical methods
(based on commutators) that are used in {\sc Section} \ref{Sec.3} to
construct algebras $\ca$ containing $\pnge$ and closed under
functional calculus.  In {\sc Section} \ref{Sec.4} we introduce the
$c_n$-calculi on manifolds with boundary and with corners. One can
embed the $c_n$-calculi, $n \ge 2$, in $\Psi^*$-algebras that consist
of smooth kernels, 
a result that is proved in {\sc Section} \ref{Sec.5}.  {\sc Section}
\ref{Sec.6} is devoted to the study of the Schwartz space $\SS(\GR)$
of a continuous family groupoid and the proof that it is closed under
functional calculus. We also define length functions for the groupoids
associated to various pseudodifferential calculi on manifolds with
corners (the $b$-calculus, the ``cusp''-calculus, or, more generally,
the $c_{n}$-calculi).

\subsubsection*{Acknowledgements:}\
The first named author is greately indebted to B. Gramsch for
introducing him to the theory of inverse closed Fr\'echet algebras and
their relations to pseudodifferential analysis.
He wants to thank R. B. Melrose and the
Massachusetts Institute of Technology where parts of this work has been
done for the invitation and  warm hospitality.

\section{Algebras closed under holomorphic functional 
calculus\label{Sec.1}}

In this section we recall some basic facts about algebras that are
closed under holomorphic functional calculus and describe a method to
generate algebras closed under holomorphic functional calculus. This
method is based on permanence properties of the closure under
holomorphic functional calculus when passing to two-sided ideals, to
quotients, or back from ideals and quotients to the algebra.

\subsection{Definitions} We begin by recalling the definition
of an algebra closed under holomorphic functional calculus.

\begin{definition}
\label{def4.1}\
Let $\cb$ be a Banach algebra with unit $e$. A not necessarily unital
subalgebra $A\subseteq\cb$ is said to be {\em closed under the
holomorphic functional calculus in $\cb$} provided for every
$a=\lambda e +x \in\cz e+A$ and all $f\in\hol{\sigma_{\cb}(a)}$, we
have $f(a)\in\cz e+A$.
\end{definition}

Here $\hol{\sigma_{\cb}(a)}$ stands for the algebra of germs of
holomorphic functions on the spectrum $\sigma_{\cb}(a)$ of $a$ with
respect to the Banach algebra $\cb$, and
\begin{equation}
\label{defhfc}
f(a):=\frac{1}{2\pi i}\int_{\gamma}f(z)(z e - a)^{-1}dz
\in\cb
\end{equation}
is the operator given by the usual holomorphic functional calculus
within the Banach algebra $\cb$. For any algebra $A$ with unit, we
shall denote by $A^{-1}$ the set of invertible elements of $A$.

\begin{remark} 
\label{rem4.2} \
The following observations are immediate. 
\begin{enumerate} 
\item Let $\cb$ be a unital $C^{*}$-algebra and 
      $A\subseteq \cb$  a symmetric subalgebra (i.e. closed under
      taking adjoints), then we have 
      $\sigma_{\cb}(a)=\sigma_{{\mathcal A}_{e}}(a)$ where
      ${\mathcal A}_{e}$ is the completion of $\cz e+A$ with 
      respect to the norm induced by $\cb$, so Definition \ref{def4.1}
        recovers  
      \cite[Definition 1, p.\ 285]{connes}. 
\item An arbitrary intersection of algebras closed under 
      holomorphic functional calculus in $\cb$ is again closed under
      holomorphic functional calculus in $\cb$.
\item If $A\subseteq\cb$ is closed under the holomorphic functional
      calculus in $\cb$, then $A$ is {\em spectrally invariant}
      in $\cb$, i.e.\ for the groups of invertible elements we have
      $$(\cz e + A)\cap\cb^{-1}=(\cz e +A)^{-1}\,,$$
      or, equivalently, 
      $j^{-1}(B^{-1})=(\cz e + A)^{-1}$ if
      $j:\cz e + A\hookrightarrow B$ denotes the natural inclusion.
\end{enumerate} 
\end{remark}

It is often useful to consider algebras closed under holomorphic
functional calculus that are in addition complete with respect to a
finer topology. For example, the class of $\Psi^{*}$-algebras
introduced by Gramsch \cite{gra84} in connection with a perturbation
theory for singular integral and pseudodifferential operators is, in
fact, the appropriate setting to describe $\cun$-phenomena of
microlocal analysis within a functional analytic framework
\cite{gra98}.

\begin{definition}
\label{def4.3}\
Let $\cb$ be a unital $C^{*}$-algebra.  A symmetric, spectrally
invariant subalgebra $\ca\subseteq\cb$, $e\in\ca$, is called a {\em
$\Psi^{*}$-algebra} in $\cb$ if, and only if, there exists a Fr\'echet
topology ${\mathcal T}_{\ca}$ on $\ca$ making the embedding $\iota:
(\ca,{\mathcal T}_{\ca})\hookrightarrow(\cb,\norm{\cdot}_{\cb})$
continuous.

In case $e\notin\ca$, the algebra $\ca$ is said to be
a {\em non-unital} $\Psi^{*}$-algebra provided
$\cz e\oplus\ca$ is a $\Psi^{*}$-algebra.
\end{definition}

Following the usual convention, in the sequel a $\Psi^{*}$-algebra
is always unital.
A unital Fr\'echet algebra $\ca$ is said to be {\em submultiplicative}
if the topology ${\mathcal T}_{\ca}$ can be generated by a countable
system $(q_{j})_{j\in\NN}$ of submultiplicative semi-norms, i.e.\
semi-norms satisfying $q_{j}(xy)\leq q_{j}(x)q_{j}(y)$ and
$q_{j}(e)=1$.

We recall a few basic facts about 
(non-unital) $\Psi^{*}$-algebra, most of them are obvious consequences
of the definition.

\begin{proposition}\label{theorem.c.f.c}\
Let $\cb$ be a unital $C^{*}$-algebra and $\ca\subseteq\cb$ 
a subalgebra.
\begin{enumerate}
\item If $\ca$ is a non-unital $\Psi^{*}$-algebra, then there
      is a Fr\'echet-topology $\mathcal{T}_{\ca}$  on $\ca$
      such that $(\ca,\mathcal{T}_{\ca})\hookrightarrow
      (\cb,\norm{\cdot}_{\cb})$ is continuous. 
\item \label{labb} 
      If $\ca$ is a $\Psi^{*}$-algebra, then the 
      group of invertible elements
       $\ca^{-1}$ is open and the inversion
      \begin{equation}\label{eq.cont.inv}
        \ca^{-1} \ni x \longmapsto x^{-1} \in \ca
      \end{equation}
      is continuous.
\item If $\ca$ is a non-unital $\Psi^{*}$-algebra, then for any 
      $a\in\ca$ there exists an analytic map
      $h:\varrho_{\cb}(a)\rightarrow\ca$ where
      $\varrho_{\cb}(a)\subseteq\cz$ is the resolvent set, such that
      $$(\lambda e - a)^{-1}= \frac{1}{\lambda}e+h(\lambda)\,.$$
\item Any $\Psi^{*}$-algebra, unital or not, is closed under functional
      calculus. If $\ca$ is non-unital and  
      $a=\lambda e +x\in\cz e \oplus \ca$ is arbitrary, then 
      $f(a)-f(\lambda)e\in\ca$ for 
      all $f \in\mathcal{O}(\sigma_{\cb}(a))$
\end{enumerate}
\end{proposition}

\begin{proof}
For \eqref{labb} it suffices to use an old result of Banach \cite{ban48} which
says that the inversion in a Fr\'echet algebra is continuous if and only if 
the group of invertible elements is a $G_{\delta}$-set.
\end{proof}

\begin{definition}
\label{Def.8}\
Let $\cb$ be a Banach algebra with unit $e$, and
$\varphi:A\rightarrow\cb$ be a morphism of algebras, which we assume
to preserve the unit if $A$ has one. Then $A$ is called {\em locally
spectral invariant with respect to $\varphi$}, if there exists
$\varepsilon>0$ such that we have
$$
           (e + \varphi(x))^{-1} \in \CC +\varphi(A)
$$
for all $x\in A$ with $\norm{\varphi(x)}_{\cb}<\varepsilon$. In that
case, we say that $A$ has {\em property} $(P_{A})$ in $\cb$, and the
morphism $\varphi$ is to be understood from the context. Moreover, $A$
is said to have {\em property} $(\wt{P}_{A})$ if
$A$ is unital and 
$\varphi^{-1}(\cb^{-1})=A^{-1}$.
\end{definition}

By \cite[Lemma 5.3]{gra84}, properties $(P_{A})$ and 
$(\wt{P}_{A})$ are closely related in many interesting
cases.

\begin{lemma}
\label{lemlocglob}
Let $\varphi:A\rightarrow \cb$ be an injective morphism of 
unital algebras. Then
\begin{enumerate}
\item If $\cb$ is a unital Banach algebra and $\varphi(A)$ is
      dense in $\cb$, then we have
      $(P_{A}) \iff (\wt{P}_{A})$.
\item If $\cb$ is  a unital $C^{*}$-algebra and $\varphi(A)$ is
      symmetric in $\cb$, then we have
      $(P_{A}) \iff (\wt{P}_{A})$.
\end{enumerate}
\end{lemma}

In other words, $(\wt{P}_A)$  is not stronger than the (apparently
weaker) condition $P_A$, provided that either $A$ is dense in $B$
or $A$ is symmetric.

We close this subsection with a lemma that helps to detect spectrally
invariant subalgebras.  It will be used in Section \ref{Sec.6} to
prove the spectral invariance of the Schwartz convolution algebra on a
continuous family groupoid.  We include here a proof that is closely
related to the one of the previous Lemma, for the sake of
completeness.

\begin{lemma}
\label{lemlaf}
Let $\cb$ be a  Banach algebra with unit $e$ and $\ca$ a 
Banach algebra such that
$e\in\ca\subseteq\cb$, $\ca$  is dense in $\cb$, and
$\varrho_{\ca}(a)=\varrho_{\cb}(a)$ 
holds for all $a\in\ca$, where
$\varrho_{T}(a)$ stands for the spectral radius of $a$ in a Banach algebra
$T$.

Then
we have

$$
        \ca\cap \cb^{-1}=\ca^{-1}\,,
$$

i.e.\ $\ca$ is spectrally invariant in $\cb$ and 
$\ca$ is closed under holomorphic functional calculus in $\cb$.
\end{lemma}

\begin{proof}
Let $a\in\ca\cap \cb^{-1}$ be arbitrary.
By the density of $\ca$ in $\cb$ there exists $y\in\ca$ with
$\norm{ay-e}_{\cb}\leq1/2$, hence
$\varrho_{\ca}(ay-e)\leq1/2$, and
$ay$ is invertible in $\ca$. 
The rest is clear. 
\end{proof}

The above statement generalizes right away to non-unital
algebras.


\subsection{Spectral invariance and ideals}

We are now going to look more closely at the question of how spectral
invariance of an algebra is related to that of its quotients, bearing
in mind Proposition \ref{theorem.c.f.c} (that a compatible Fr\'echet
topology on a spectrally invariant subalgebra implies stability under
holomorphic functional calculus). We start by describing the general
setting.

{F}rom now on and throughout this section, $\cb$ will be a Banach algebra
with unit $e$, $A\subseteq \cb$ will be a subalgebra with $e\in A$,
$\cj\subseteq\cb$ will be a proper, closed two-sided ideal in $\cb$,
and $I\subseteq A$ will be a two-sided ideal in $A$ with $I\subseteq
\cj$. Then the map
\begin{equation}\label{eq.phi}
        \varphi:A/I\longrightarrow\cb/\cj:a+I\longmapsto a+\cj
\end{equation}
is a well-defined homomorphism of unital algebras.  Also, note that
$\varphi$ is one-to-one if, and only if, $\cj\cap A=I$.  Thus, we have
a commutative diagram with exact rows
\begin{equation}
\label{diag1}
\begin{array}{ccccccccc}
0 & \longrightarrow & I & \longrightarrow & A & \longrightarrow & A/I
& \longrightarrow & 0 \\ && \Big\downarrow &&\Big\downarrow &&
{}\hspace{1ex}\Big\downarrow\varphi\\ 0 & \longrightarrow & \cj &
\longrightarrow & \cb & \longrightarrow & \cb/\cj & \longrightarrow &
0 \\
\end{array}\,.
\end{equation}

The next theorem relates the different properties $(P_I), (P_A),
(P_{A/I})$, and $(\wt{P}_{A/I})$, to one another.  A special case of
part (a) can be found in \cite[Appendix]{bues}, (b) is from
\cite[Section 5]{gra90}.

\begin{theorem}\label{thm4.02}\
Let $A,I, \cb$, $\cj$, and $\varphi$ be as above (Equation
\eqref{eq.phi}). We write $(\wt{P})$ instead of $(\wt{P}_{A/I})$, for
simplicity; the property $P_{I}$ is considered with respect to the
obvious morphism $I\rightarrow\cz e\oplus \cj$. Then we have
\begin{enumerate} 
\item If $I\subseteq\cj$ is dense, then $(P_I)$ together with
      $(\wt{P})$ imply $(P_A)$.
\item If $I\subseteq\cj$ is dense, then $(P_A)$ implies $(P_{A/I})$.
\item If $I\subseteq\cj$ and $A\subseteq\cb$ are dense, then
      $(P_A)$ implies $(\wt{P})$.
\item $(P_A)$ implies $(P_I)$.
\item $(\wt{P})$ implies $(P_{A/I})$; moreover, if 
      $\varphi:A/I\longrightarrow\cb/\cj$ is one-to-one
      and $\varphi(A/I)$ is dense in $\cb/\cj$, then we have
      $(P_{A/I})$ if and only if we have $(\wt{P})$.
\end{enumerate} 
\end{theorem}

\begin{proof}\ Let $\varepsilon >0$ always be the constant from Definition
\ref{Def.8}.

(a) By the continuity of the inversion in the unital Banach algebra
$\cb/\cj$, we can find $0<\delta<\varepsilon/2<1$ such that
$\norm{(b+\cj)^{-1}-(e+\cj)}_{\cb/\cj}<\varepsilon/8$,
for all $b \in \cb$ satisfying
$\norm{(b+\cj)-(e+\cj)}_{\cb/\cj}<\delta$.
      
Consider now $a\in A$ with $\norm{a-e}_{\cb}<\delta$. From
$\norm{(a-e)+\cj}_{\cb/\cj}<\delta$ and $(\wt{P})$, we obtain that
there exists $a_{1}\in A$ with $aa_{1}-e=:x_{1}\in I$ and
$\norm{a_{1}-e+\cj}_{\cb/\cj}<\varepsilon/8$. Using the density of $I$
in $\cj$, we find $x_{2}\in I$ with
$\norm{a_{1}-e+x_{2}}_{\cb}<\varepsilon/4$, hence
$$
        \norm{x_{1}+ax_{2}}_{\cb}\leq\norm{a}_{\cb}
      \norm{a_{1}-e+x_{2}}_{\cb}+\norm{a-e}_{\cb}<
      \varepsilon\,,
$$ 
and by $(P_I)$ we get $z\in I$ with $(e+x_{1}+ax_{2})^{-1}=e+z$, i.e.\ 
$$
        a(a_{1}+x_{2})(e+z)=e,
$$ 
which gives $a\in A^{-1}$, and completes the proof of (a).
     
(b) Let $0<\delta<\varepsilon/2$, and $a+ I\in A/I$ be with
$\norm{\varphi(a+I)-(e+\cj)}_{\cb/\cj}<\delta$. By the density of $I$
in $\cj$, there exists $x\in I$ with $\norm{a+x-e}_{\cb}<\varepsilon$,
hence $(P_A)$ gives $a_{1}\in A$ with $(a+x)a_{1}=e=a_{1}(a+x)$. We
thus obtain 
$$
        \varphi(a+I)^{-1}=\varphi(a_{1}+I)\in\varphi\left(A/I\right),
$$
and hence (b) is proved.

(c) Let $a+I\in A/I$ be with
$\varphi(a+I)\in\left(\cb/\cj\right)^{-1}$. Thus, there exists
$b\in\cb$ and $y\in\cj$ with $ab-y-e=0$. By the density assumption we
obtain $x\in I$ and $a_{1}\in A$ with
$\norm{aa_{1}-x-e}_{\cb}<\varepsilon$, hence $aa_{1}-x\in A^{-1}$ by
$(P_A)$. Let $a_{2}\in A$ be such that
$e=(aa_{1}-x)a_{2}=aa_{1}a_{2}-xa_{2}$.  Consequently,
$a_{1}a_{2}+I\in A/I$ is a right-inverse of $a+I$ in $A/I$. Similarly,
we obtain also a left-inverse of $a + I$, which gives
$a+I\in\left(A/I\right)^{-1}$, and completes the proof of (c).

To prove (d), let $\varepsilon>0$ and $x\in I$ be with
$\norm{x}_{\cb}<\varepsilon$. Then $e + x$ is invertible in $A$
because of $(P_A)$. Then, $(P_I)$ is a consequence of $(P_A)$ and the
identity
$$
        (e+x)^{-1}-e=-x+\underbrace{(e+x)^{-1}}_{\in A} x^{2}\in I\,.
$$ 

Finally, (e) is a simple, straightforward computation using Lemma
\ref{lemlocglob}.
\end{proof}

We are mostly interested in the question when
$A$ is spectrally invariant in $\cb$. The following
special case of Theorem \ref{thm4.02} will be used in the sequel.

\begin{corollary}
\label{coreins}\
Let $\cb$ be a unital $C^{*}$-algebra, and $\cj$, $A$, and $I$ as above,
but additionally symmetric with respect to the $*$-operation in
$\cb$. Assume that $I$ is dense in $\cj$. Then $A$ is spectrally invariant
in $\cb$ provided
\begin{enumerate}
\item $\cz e\oplus I$ is spectrally invariant in $\cz e\oplus \cj$, and
\item  $\varphi^{-1}\left((B/\cj)^{-1}\right)=\left(A/I\right)^{-1}$.
\end{enumerate}
\end{corollary}

\section{Semi-ideals\label{Sec.2}}

We now turn to a 
method of constructing algebras closed under
holomorphic functional calculus, or, more generally,
$\Psi^*$-algebras.  More precisely, in order to be able to deal with
non-unital algebras we need to study 
composition with possibly
unbounded operators and the semi-ideals generated by 
this process.

\subsection{Definitions} We now introduce semi-ideals.
\begin{definition}
\label{def4.4}\
A subspace $J\subseteq B$ of a unital algebra $B$ is said to be a {\em
semi-ideal} in $B$ provided that we have $xby\in J$ for all $x,y\in J$
and all $b\in B$.
\end{definition}

\begin{remark}
\label{rem4.5}
\quad
\begin{enumerate}
\item A left, right or two-sided ideal is obviously a semi-ideal.
\item Let $B$ be a unital algebra and $J\subsetneq B$ be a proper
      semi-ideal. Then we have $e\notin J$ and
      $J\cap B^{-1}=\emptyset$
      because otherwise we would contradict the identities
       $b=ebe$ and $e=x(x^{-1})^{2}x$.
\end{enumerate}
\end{remark}

\begin{proposition}\label{prop4.6}\
Let $\ca\subseteq\cb$ be a $\Psi^{*}$-algebra, and $J\subsetneq\ca$ be 
a proper semi-ideal.
\begin{enumerate}
\item $0\in\sigma_{\cb}(x)$ for all $x\in J$.
\item $J$ is closed under the holomorphic functional calculus in $\cb$.
\end{enumerate}
\end{proposition}

\begin{proof}\
The first statement is an immediate consequence of Remark
\ref{rem4.5}(b) and the spectral invariance of $\ca$ in $\cb$.  If
$a=\lambda e+x\in\cz e\oplus J$ and $f\in\hol{\sigma_{\cb}(a)}$ are
arbitrary, then we have $\lambda\in\sigma_{\cb}(a)$ by (1), 
so for all $\mu\notin\sigma_{\cb}(a)$ we have 
$$
        (\mu e-a)^{-1} = \frac{1}{\mu - \lambda}e +
        \frac{1}{(\mu-\lambda)^{2}} x + x\frac{1}{(\mu -
        \lambda)^{2}}(\mu e - a)^{-1}x
$$
together with the definition (\ref{defhfc}) of $f(a)$ lead to
$$
        f(a)= f(\lambda)e+f'(\lambda)x + x\, 
               \underbrace{\frac{1}{2\pi i}\int_{\gamma}
        \frac{f(\mu)}{(\mu-\lambda)^{2}} (\mu e -a )^{-1}
        d\mu}_{\in\ca} \, x,
$$
which gives $f(a)-f(\lambda)e\in J$, and completes the proof.
\end{proof}

\begin{remark}
\label{rem4.7}\
A careful inspection of the above proof shows that (a) holds for any
spectrally invariant, unital subalgebra $A\subseteq \cb$, and (b) remains
true for any spectrally invariant, continuously embedded, sequentially 
complete, locally convex topological algebra $A$ with continuous inversion 
and $e\in A$. 
\end{remark}

If $\cb$ is in addition a $C^{*}$-algebra, we can assume
without loss of generality that algebras
closed under holomorphic functional calculus are symmetric. More precisely,
we have the following result.

\begin{lemma}\label{lem4.8}\
If $A\subseteq \cb$ is closed under the holomorphic functional
calculus in 
the $C^{*}$-algebra $\cb$, then $A_{*}:=\{a\in A: a^{*}\in A\}$ is symmetric
and closed under the holomorphic functional calculus in $\cb$.
\end{lemma}

\begin{proof}\
It is sufficient to note that we have
$f(a)^{*}=f^{*}(a^{*})$ if $f\in\hol{\sigma_{\cb}(a)}$ and
$[f^{*}:z\longmapsto\overline{f(\bar{z})}]\in\hol{\sigma_{\cb}(a^{*})}$.
\end{proof}

\subsection{A commutator method}
We shall use several procedures to construct subalgebras closed under
holomorphic functional calculus. The first one leads to
$\Psi^{*}$-algebras using commutator methods, whereas the second one
produces a semi-ideal, hence also an algebra closed under functional
calculus. We begin by recalling the construction of submultiplicative
$\Psi^{*}$-algebras using commutators with closed, symmetric
operators.  These techniques were first used in \cite{beals77} for the
characterization of pseudodifferential operators on $\RR^{n}$, and
later on investigated systematically in \cite{guw}, for instance. Here
we follow the presentations in \cite{guw,aa}.  Let us start with a
description of the general setting.

Let ${\mathcal K}$ be a Hilbert space,
$(\ca,(\norm{\cdot}_{j})_{j\in\NN})$ be a submultiplicative
$\Psi^{*}$-algebra in $\lck$ with $\norm{\cdot}_{1} =
\norm{\cdot}_{\lck}$. Also, let $\ft$ be a finite set of densely
defined, closed, symmetric operators $T:{\mathcal
K}\supseteq\de{T}\longrightarrow {\mathcal K}$.  By \cite[Lemma
2.17]{aa} each $T\in\ft$ induces a closed ${}^{*}$-derivation
$$
        \dt:\ca\supseteq\de{\dt}\longrightarrow\ca.
$$
Here, we have $a\in\de{\dt}$ if, and only if, $a(\de{T})\subseteq
\de{T}$, there exists $\dt(a)\in\ca$ with
$\dt(a)\varphi=i\left(Ta\varphi - aT\varphi\right)$, for all
$\varphi\in\de{T}$, and the same is also true for $a^{*}$.

Furthermore, as shown in \cite[Section 2.2]{aa}, the finite set
$\{\dt:T\in\ft\}$ leads to the following  
scale of symmetric subalgebras of $\ca$
\begin{eqnarray*}
\pft{0} & := & \ca\,, \\
\pft{1} & := & \bigcap_{T\in\ft}\de{\dt}\,,\\
\pft{r} & := & \{a\in\pft{r-1}:\dt(a)\in\pft{r-1}\mbox{ for all }
                    T\in\ft\}, \, r\geq 2\,.
\end{eqnarray*}
Each of the algebras $\pft{r}$ is endowed with a system of semi-norms,
namely, 
\begin{eqnarray*}
q_{0,j}(a) & := & \norm{a}_{j}\,,
                  \mbox{ for } a\in\pft{0}, j\in\NN, \mbox{, and }\\
        q_{r,j}(a)&:=&  q_{r-1,j}(a)+\sum_{T\in\ft}
        q_{r-1,j}(\dt(a))\,,\mbox{ for }
        a\in\pft{r}\,, r,j\in\NN\,.
\end{eqnarray*}
Finally, we endow the algebra $\puft:=\bigcap_{r=0}^{\infty}\pft{r}$ with the 
system of semi-norms $(q_{r,j})_{r,j\in\NN}$.

Similarly, the set $\ft$ induces the scale of  so-called 
{\em $\ft$-Sobolev spaces} by
\begin{eqnarray*}
\hft{0} & :=  & {\mathcal K}\\
\hft{1} & := & \bigcap_{T\in\ft}\de{T}\\
\hft{r} & := & \{x\in \hft{r-1}Tx\in\hft{r-1}\mbox{ for all }
                  T\in\ft\}\,, r\geq2\,.
\end{eqnarray*}
As above, the spaces $\hft{r}$ are endowed with the iterated graph
norms with respect to the system $\ft$, i.e.\ $p_{0}(x) :=
\norm{x}_{{\mathcal K}}$ and
$$
        p_{r}(x):=p_{r-1}(x)+\sum_{T\in\ft}p_{r-1}(Tx)\,,
        x\in\hft{r}\,, r\geq1\,.
$$
The intersection $\huft:=\bigcap_{r=0}^{\infty}\hft{r}$ is endowed
with the system of norms $(p_{r})_{r\in\NN}$.

The main properties of this construction are summarized in the next
theorem. For a proof see \cite[Theorem 2.24]{aa}, or 
\cite[Section 2]{guw} for the  special case $\ca=\lck$.

\begin{theorem}\label{th4.10}\
The algebra $(\puft,(q_{r,j}))$ is a submultiplicative
$\Psi^{*}$-algebra in $\lck$ with $\puft\subseteq\ca$. The
$\ft$-Sobolev-spaces $\hft{r}$ are Hilbert spaces, $\huft$ is a
projective limit of a sequence of Hilbert spaces, and for each
$r\in\NN\cup\{\infty\}$ the natural map
$$
        \pft{r}\times\hft{r}\longrightarrow
        \hft{r}:(a,\varphi)\longmapsto a(\varphi)
$$ 
is bilinear and continuous.
\end{theorem}

\subsection{Commutators and semi-ideals}
The second construction that we shall need associates to 
the algebra $\ca$ and the system $\ft$ a semi-ideal in $\lck$.

Indeed, let $\cjt{0}:=\ca$, and denote by $\cjt{1}$ the space of all
$x\in\cjt{0}$ such that, for all $T,T_{1},T_{2}\in\ft$, we have
\begin{enumerate}
\item $x({\mathcal K})\subseteq\de{T}$ and $\olt(x):=Tx\in\cjt{0}$.
 \item There exists $\ort(x)\in\cjt{0}$ with
      $\ort(x)f=xTf$ for all $f\in\de{T}$.
\item $x({\mathcal K})\subseteq\de{T_{1}}$, and there is
$\olrt(x)\in\cjt{0}$ with $\olrt(x)f=T_{1}xT_{2}f$ for all
$f\in\de{T_{2}}$.
\end{enumerate}
Moreover, let $\cjt{k+1}$ be the space of all $x\in\cjt{k}$ such
that we have 
$\olt(x)$, $\ort(x)$, $\olrt(x)\in\cjt{k}$ for all
$T,T_{1},T_{2}\in\ft$.

We endow the spaces $\cjt{k}$ with the following systems 
$(p_{j,k})_{j\in\nz_{0}}$ of norms: let $p_{j,0}(x)  : = 
\norm{x}_{j}$, for $x\in\cjt{0}$, and let
$$
        p_{j,k+1}(x)  : = 
        p_{j,k}(x)+\sum_{T\in\ft}\left(p_{j,k}(\olt(x))
        +p_{j,k}(\ort(x))\right) 
        + \sum_{T_{1},T_{2}\in\ft}p_{j,k}(\olrt(x))\,,
$$
for $x\in\cjt{k+1}$. Moreover, 
the projective limit $\cjtu:=\bigcap_{k=0}^{\infty}\cjt{k}$ is endowed 
with
the projective topology given by the system of norms
$(p_{j,k})_{j,k\in\nz_{0}}$.

Let us collect the main properties of this construction in the
following theorem.

\begin{theorem}\label{thm4.12}\
Let $(\ca,(\norm{\cdot}_{j})_{j\in\NN})$ be a submultiplicative
$\Psi^{*}$-algebra. The we have for $k\in\nz_{0}\cup\{\infty\}$:
\begin{enumerate}
\item $\left(\cjt{k},(p_{j,k})\right)$ is a submultiplicative
        Fr\'echet algebra. The canonical embedding
        $\cjt{k}\hookrightarrow\cjt{0}$ is continuous.
\item $\puft\cjt{k}\subseteq\cjt{k}$ and
      $\cjt{k}\puft\subseteq\cjt{k}$;
      the two canonical bilinear maps,
$$
        \puft\times\cjt{k}\longrightarrow\cjt{k}
      \mbox{ and }
      \cjt{k}\times\puft\longrightarrow\cjt{k},
$$
      are jointly continuous.
\item $\cjt{k}$ is a semi-ideal in the $\Psi^{*}$-algebra
      $\ca$; in particular,
      $\cjt{k}$ is closed under the holomorphic functional calculus in 
      $\lh$. Moreover, the canonical map
      \begin{equation}
      \label{maneq1}
      \cjt{k}\times\ca\times\cjt{k}\longrightarrow\cjt{k}
      \end{equation}
      is jointly continuous.
\end{enumerate}
\end{theorem}

\begin{proof}\
The proofs are by induction with respect to $k$.
Since the arguments for the steps from $k$ to $k+1$ are the same as
for $k=1$ the proofs of the steps are omitted.

      For (a), let us first 
      assume that we have $x,y\in\cjt{1}$. Then we have 
      $xy\in\cjt{1}$ with
      $$
      \olt(xy) = \olt(x)y\,,\, \ort(xy)=x\ort(y)\mbox{, and }
      \olrt(xy)=\olta{T_{1}}(x)\orta{T_{2}}(y)$$ for all
      $T,T_{1},T_{2}\in\ft$; the submultiplicativity is now immediate
      whereas for the completeness of $\cjt{1}$ we have to use the
      closedness of the operators $T\in\ft$.

      For (b),  note that for $a\in\puft$ and $x\in\cjt{1}$ 
      a straightforward computation gives  
      $ax\in\cjt{1}$ and $xa\in\cjt{1}$ with
      \begin{eqnarray*}
      \olt(ax) & = & a\olt(x)-i\dt(a)x \\ 
      \ort(ax) & = & a\ort(x)\\ 
      \olrt(ax) & = & a\olrt(x)-i\delta_{T_{1}}(a)\orta{T_{2}}(x) \\ 
      \olt(xa) & = & \olt(x)a\\
      \ort(xa) & = & \ort(x)a-ix\dt(a)\\
      \olrt(xa) & = & \olrt(x)a+i\olta{T_{1}}(x)\delta_{T_{2}}(a)
      \end{eqnarray*}
      for all $T,T_{1},T_{2}\in\ft$.

     For the semi-ideal property (c) it suffices to note that  
      for $a\in\ca$ and $x,y\in\cjt{k}$ we have
      $$\olt(xay)=\olt(x)ay\,,\,
        \ort(xay)=xa\ort(y) \mbox{, and }
      \olrt(xay)=\olta{T_{1}}(x)a\orta{T_{2}}(xy)$$
      for all $T,T_{1},T_{2}\in\ft$, which gives the joint continuity
      of (\ref{maneq1}) as well.
\end{proof}

\begin{remark}
\label{rem4.13'}
Note that it is not clear, and in general not true that the spaces
$\cjtk$ are symmetric subspaces of $\lck$. However, we easily obtain
this property by considering the spaces
\begin{equation}
\label{symsemi}
\cjtk_{*}:=\{x\in\cjtk:x^{*}\in\cjtk\}\,.
\end{equation}
It is straightforward to check that Theorem \ref{thm4.12} remains true
also for the smaller spaces $\cjtk_{*}$. By slight abuse of notation
we will 
sometimes write $\cjtk$ for the spaces $\cjtk_{*}$.
\end{remark}

\begin{corollary}\label{cor.psi.J}\
We have that $\cjt{m}$ is a non-unital $\Psi^*$-algebra, for any $m
\in \NN \cup \{\infty\}$.
\end{corollary}

\section{An extended pseudodifferential calculus}
\label{Sec.3}

Recall that the notion of a \cfg, defined in \cite{Paterson}
generalizes that of $\cinfo$-foliations as considered in
\cite{connesF}.  More precisely, a {\em continuous family groupoid} is a
locally compact topological groupoid such that $\GR$ is covered by
some open subsets $ \Omega$ such that the following conditions are satisfied:
\begin{itemize}
\item each chart $\Omega$ is homeomorphic to two open subsets of
  $\R^k\times \GR^{(0)} $,
  $ T_d\times U_d$ and $ T_r\times U_r$
   such that the following diagram is commutative:

\begin{equation}
\label{coord}
\xymatrix{ & T_r\times U_r \ar[dl]& \Omega \ar[l]_<<<<{\simeq}
        \ar[dl]^r \ar[r]^>>>>>{\simeq} \ar[dr]_d& T_d\times U_d
        \ar[dr] \\ U_r & r(\Omega)\ar[l]_= && d(\Omega) \ar[r]^= &
        U_d}
\end{equation}

\item  each coordinate change (for $r$ resp.\  $d$)
      is given by $(t,u) \mapsto (\phi(t,u),u)$ where
  $\phi$ is of class $\cinfo$, i.e.\  $u \mapsto \phi(.,u)$ is a
  continuous map from $U_{*}$ 
  to ${\mathcal C}^{\infty}(T_{*},T_{*}')$, $*=d,r$.
\end{itemize}
In addition, one requires that the composition and the inversion be $\cinfo$ 
morphisms. For the sake of simplicity, we will always assume
that the space $M:=\GR^{(0)}$ of units is compact.

\subsection{Groupoid algebras}
To any \cfg\ $\GR$,  there is associated
 an algebra of \pd\ operators:
namely, let $\Psi^{m,0}(\GR)$ be the space of continuous, 
uniformly supported, 
invariant families of \pd\ operators $(P_{x})_{x\in M}$
on the fibers of the \gr, see \cite{lmn}.
For later purposes, note that this implies in
particular the following:
Let $\Omega\cong T\times U$ be a chart as in \eqref{coord} and
$\varphi\in\ccinfo(\Omega)$. Then there exists
$p\in\mathcal{C}_{c}(U,S^{m}(T;\rz^{k}_{\eta}))$ such that
for each $x\in U$, 
$\varphi P_{x}\varphi$ corresponds to
the pseudodifferential operator
$p(x,y,D_{y})$ on $T$. 

In this section, we are going to show how the methods of the previous
section can be used to construct algebras closed under holomorphic
functional calculus that contain the algebras $\Psi^{-\infty,0}(\GR)$
and $\Psi^{0,0}(\GR)$ of pseudodifferential operators on the groupoid
and share some of their algebraic and analytic properties.  As it
might be expected, it is difficult to construct, in general, algebras
closed under holomorphic functional calculus that retain all geometric
properties of the given groupoid.  Nevertheless, we do construct
algebras closed under holomorphic functional calculus that retain at
least some of the geometric properties of the groupoid.

First, let us fix some notations. Throughout this section, 
${\mathcal K}_{x}:=\lzgx$  stands for the completion of the space 
$\cunc(\gx;\rsed)$ with respect to the sesquilinear pairing
$$(f,g)_{x}:=\int_{\gx}f(\gamma)\overline{g(\gamma)}
\,.$$
Moreover, the Hilbert spaces ${\mathcal K}_{x}$, $x\in M$, can  be glued
together to a new, big Hilbert space $\cH$ 
containing all necessary
information for us. 

To be more precise, fix a positive density
$\nu\in\cun( M,\Omega)$  on $ M$, and denote the induced measure
again  by $\nu$. 
Note that each $f\in\ccinfo(\GR;\rsed)$ induces a section
$$
        \tilde{f}: M\longrightarrow\prod_{x\in M}{\mathcal
        K}_{x}:x\longmapsto f|_{\gx}\in\cunc(\gx;\rsed)\subseteq {\mathcal
        K}_{x}\,,
$$ 
where, as usual, ``section'' simply  means $\tilde{f}(x)\in
{\mathcal K}_{x}$ for all $x\in M$.  Now, for any
$f,g\in\ccinfo(\GR;\rsed)$, the function 
$ 
        M\ni x \mapsto (\tilde{f}(x),\tilde{g}(x))_{x}
$ 
is continuous, hence $\nu$-measurable.  Consequently, the set $\fm$ of
all sections $h: M\rightarrow\prod_{x\in M}{\mathcal K}_{x}$ such
that the map $M\ni x\mapsto (\tilde{f}(x),h(x))_{x}$ is $\nu$-measurable
for all $f\in\ccinfo(\GR,\rsed)$ induces the structure of a
$\nu$-measurable field of Hilbert spaces on the family $({\mathcal
K}_{x})_{x\in M}$ \cite[Definition 2.1.3.1, Proposition
2.1.4.4]{dixn}. The set $\cH$ of all $h\in\fm$ satisfying 
$\int_{M}\norm{h(x)}^{2}_{{\mathcal K}_{x}} d\nu(x)<\infty$ is in fact a
Hilbert space and we  write $\cH=:\dirint_{ M}{\mathcal K}_{x}d
\nu(x)\,.$ The space $\ccinfo(\GR;\rsed)$ is then a dense subspace of $\cH$.

Recall that an operator $P\in\lch$ is said to be {\em decomposable}
provided there exists a family $(\wh{P}(x))_{x\in M}$ of operators
$\wh{P}(x)\in\lhx$ such that, for any $h\in\fm$  
$$
        (Ph)(x)=\wh{P}(x)h(x),
$$ 
for all $x\in  M$, and 
$[x\mapsto \norm{\wh{P}(x)}_{\lckx}]\in L^{\infty}( M;\nu)$
\cite[Definition 2.2.3.2]{dixn}.
As usual, we write in that case
$P=\dirint_{ M}\wh{P}(x)d\nu(x)$.
A straightforward computation gives \cite[Proposition 2.2.3.2]{dixn}
$$
        \norm{P}_{\lch}=
\operatorname{ess-sup}_x\norm{\wh{P}(x)}_{\lckx}.
$$
The set of all decomposable operators is in fact  a
$C^{*}$-subalgebra, which we denote by $\ldch\subseteq\lch$.

Finally, for any $\gamma\in\GR$,
the operators
$$
        U_{\gamma}:\cunc(\GR_{d(\gamma)};\rsed)
        \longrightarrow\cunc(\GR_{r(\gamma)};\rsed), \quad
        U_{\gamma}(f) (\gamma') = f (\gamma'\gamma),     
$$
extend by continuity to isometric isomorphisms $U_{\gamma}:{\mathcal
K}_{d(\gamma)}\longrightarrow {\mathcal K}_{r(\gamma)}$, and hence
induce an action of the groupoid $\GR$ on $\cH$.  Decomposable
operators $P\in\ldch$ with
$\wh{P}(r(\gamma))U_{\gamma}=U_{\gamma}\wh{P}(d(\gamma))$ for all
$\gamma\in\GR$ are called {\em invariant} with respect to the action
of the groupoid, and we denote the $C^{*}$-subalgebra of all invariant
operators by $\ldgch$.

By the results of \cite[Section 3]{lmn}, the regular representations
$\pi_{x}$, $x\in  M$,  of the algebra $\pnge$ fit together to a faithful
${}^{*}$-representation 
$$
        \pi_{r}:\pnge\longrightarrow\ldgch: 
        \pi_r(P) = \int_M^\oplus \pi_x(P)\,d\nu(x)
$$
of unital algebras.
In the sequel, we will identify $\pnge$ with its image under
$\pi_{r}$ and construct subalgebras $A\subseteq\lch$
closed under holomorphic
functional calculus in $\lch$ and containing  $\pnge$ as a subalgebra. 
Because of Lemma \ref{lem4.8} and the following Lemma, we can always
assume that the algebras $A$ are symmetric subalgebras of $\ldgch$.

\begin{lemma}\label{lem5.1}\
Let $A\subseteq\ldch$ be closed under  holomorphic functional
calculus in $\lch$. Then $A^{\GR}:=A\cap\ldgch$ is also
$\GR$-invariant and closed under  holomorphic functional calculus
in $\lch$.
\end{lemma}

\begin{proof}\
This follows because $\ldgch$ is a $C^{*}$-subalgebra of $\lch$.
\end{proof}

As a first step towards constructing algebras
closed under holomorphic functional calculus,
we are going to show that we can reduce the problem of 
finding such algebras $A$ essentially to the construction of algebras
that contain $\pmunge$ and are closed under the holomorphic functional 
calculus in $\lch$. Indeed, let 
$$
        \cj:=\rideal{\GR}=\overline{\pmunge}^{\lch},
$$
$$
        \cb:=\ralg{\GR}=\overline{\pnge}^{\lch},
$$
and suppose that we have a
subspace $I=I^*\subseteq\ldgch$ with the following properties
\begin{eqnarray}
\label{con1} &&
\pmunge\subseteq I 
\subseteq \cj
 \,,\\
\label{con2}
&&I \mbox{ is a } \pnge \mbox{-left and -right module, and } \\
\label{con3}
&& I \mbox{ has property } (P_I) \mbox{ in } \cb\,.
\end{eqnarray}

\begin{theorem}\label{propred}\
Let $I$ be as above and
$A:=\pnge +I\subseteq \ldgch$. Then we have
$$(\cz id_{\cH}+A)\cap\lch^{-1}=(\cz id_{\cH}+A)^{-1}\,.$$
In particular, if there exists a Fr\'echet topology on $I$ making the
$\pnge$-module action as well as the embedding $I\hookrightarrow\lch$
continuous, then $A$ is closed under the holomorphic functional
calculus in $\lch$ and 
is a $\Psi^{*}$-algebra containing $\pnge$.
\end{theorem}

\begin{proof}\
Because of (\ref{con2}), the space $I\subseteq A$ is an ideal in the
algebra $A$. By Lemma \ref{lemlocglob} and Theorem \ref{thm4.02}(a), it
suffices to prove that $A/I$ has property $(\wt{P}_{A/I})$ in
$\cb/\cj$. So, let $\bar{a}=a+I\in A/I$ be such that
$\varphi(\bar{a})\in(\cb/\cj)^{-1}$.  Without loss of generality, we
may assume $a\in\pnge$. On the other hand, the homogeneous principal
symbol map $\sigma_{0}$ induces an isomorphism
$\wh{\sigma}_{0}:\cb/\cj\longrightarrow\mathcal{C}(S^{*}(\GR))$, with
$\sigma_{0}(a)\in \ccinfo(S^{*}(\GR))^{-1}$.  Using the exactness of
the sequence
\begin{equation}
\label{syse1}
0 \lra \Psi^{m-1,0}(\GR;E) \lra \Psi^{m,0}(\GR;E) 
\stackrel{\sigma_m}{-\!\!\!-\!\!\!\lra}
\II \lra 0
\end{equation}
and  the asymptotic completeness of
$\pnge$, we obtain $b\in\pnge$ with
$$
        id_{\cH}-ab\in\pmunge\,,\;id_{\cH}-ba\in\pmunge.
$$
Thus, $\bar{a}$ is invertible in $\pnge/\pmunge$, hence also
in $\pnge/(I\cap\pnge)$ because of $\pmunge\subseteq I\cap\pnge$.
This gives the property $(\wt{P}_{A/I})$, and completes the proof.
\end{proof}

The following Lemma enables us to apply the methods developed in the
previous section to pseudodifferential operators on the groupoid
$\GR$.

\begin{lemma}\label{lemmamin}\
Let $T\in\pmge$ be arbitrary. Then the unbounded operator
$$
        T:\cH\supseteq\CIc(\GR;\rsed)\longrightarrow\cH
$$ 
is closable.  Moreover, if $T=T^*\in\pmge$, then the closure of $T$ is
symmetric.
\end{lemma}

For notational simplicity, we shall not distinguish between
$T\in\pmge$ and its minimal closed extension 
$\overline{T}:\cH\supseteq\de{\overline{T}}\rightarrow\cH$ in the
sequel.

\begin{proof}\
Let $(f_{j})_{j\in\NN}$ be a sequence in $\CIc(\GR;\rsed)$ with
$f_{j}\longrightarrow 0$ in $\cH$ and $Tf_{j}\longrightarrow f$ in
$\cH$ for some $f\in\cH$. We have to show $f=0$.  Let $T^{*}\in\pmge$
be the formal adjoint of $T$. Then we
obtain, for all $h\in\ccinfo(\GR,\rsed)$
\begin{eqnarray*}
        <f,h>_{\cH} & = & \lim <Tf_{j},h>_{\cH}\\ &= & \lim
        \int_{M}<T_{x}f_{j}|_{\GR_{x}},h|_{\GR_{x}}>_{{\mathcal
        K}_{x}} d\nu(x)\\ & = & \lim
        \int_{M}<f_{j}|_{\GR_{x}},T^{*}_{x}h|_{\GR_{x}}>_{{\mathcal K}
        _{x}}d\nu(x)\\ &=& \lim <f_{j},T^{*}h>_{\cH}= 0\,.
\end{eqnarray*}
Since $\ccinfo(\GR,\rsed)$ is dense in $\cH$, this gives $f=0$. The
rest is clear.
\end{proof}

Let us briefly outline how the operator theoretic methods from the
previous sections are used to embed the algebra $\pnge$ of
pseudodifferential operators of order $0$ in an algebra that is closed
under holomorphic functional calculus.

We start with the construction of a submultiplicative
$\Psi^{*}$-algebra $\can\subset\lch$ with $\pnge\subseteq\can$ by
using a variant of the commutator methods described in Theorem
\ref{th4.10}. Without loss of generality, we can further assume that
$\can\subseteq\ldgch$ and $\pnge$ dense in $\can$
\cite[Corollary 2.5]{aa}.  Then, we take a finite set
$\ft\subseteq\pege$ of formally self-adjoint first order
pseudodifferential operators and construct the semi-ideal $\cjtu$ and
the $\Psi^{*}$-algebra $\puft$ as described in Theorem
\ref{thm4.12}. By Theorem 1 of \cite{lmn} we have $\pmunge \subseteq
\cjtu$ and $\pnge \subseteq \puft$, hence (\ref{con1}), (\ref{con2}),
and (\ref{con3}) hold for the symmetrized semi-ideal $I:=\cjtu_{*}$ by
Theorem \ref{thm4.12}, thus $\cae:=\pnge+\cjtu_{*}$ is a
(submultiplicative) $\Psi^{*}$-algebra containing $\pnge$ by Theorem
\ref{propred}. Again by \cite[Corollary 2.5]{aa}, we can even assume
that $\pnge$ is dense in $\cae$.

Since the arbitrary intersection of algebras closed under
functional calculus in $\lch$ is still closed under  holomorphic
functional calculus in $\lch$, we can admit also countable subsets
$\ft\subseteq\pege$, and then we obtain, by taking the intersection
over all possible choices involved in the construction of the algebra
$\cae$, an algebra $\caun$ that is independent of any choices, is
closed under  holomorphic functional calculus in $\lch$, and
contains $\pnge$. It remains to describe some properties of the
elements of the algebra $\caun$.

\subsection{Properties of the algebra $\can$} 
Let $\chi:\Omega\stackrel{\cong}{\longrightarrow}
U\times T\subseteq M \times\rz^{k}_{y}$ be a 
chart of $\GR$ as in \eqref{coord}, and 
$\varphi_{0},\psi_{0}\in\ccinfo(\Omega)$ be arbitrary. We
define the {\em local symbol} of $a\in\ldch$ on $\Omega$ with respect to
$\varphi_{0}$ and $\psi_{0}$ by
$$
        \sva(x,y,\eta):= e^{-iy\eta}\psi_{0}(x,y)
        \wh{a}(x)[y'\longmapsto\varphi_{0}(x,y')e^{iy'\eta}](y)
$$
for $x\in U$, $y\in T$, and $\eta\in\rz^{k}_{\eta}$. Thus, we get
\begin{eqnarray}\label{l2est}
        \sup_{x,\eta}\int_{T}|\sva(x,y,\eta)|^{2}dy & \leq &
        c(\varphi_{0},\psi_{0})\sup_{x}\norm{\wh{a}(x)}^{2}_{\lckx}\\
        \nonumber &=&c(\varphi_{0},\psi_{0})\norm{a}^{2}_{\ldch}\,.
\end{eqnarray}
Since the operators $i\partial_{y_{s}}$, $iy_{r}\partial_{y_{s}}$,
$r\neq s$, $i y_{s}\partial y_{s}+\frac{i}{2}$, and $M_{y_{r}}$ after
multiplication with some cut-off function can be realized as closed
symmetric operators on $\cH$, we can consider the corresponding
submultiplicative $\Psi^{*}$-algebra $\cavs$ with respect to the
$C^{*}$-algebra $\ldch$ as in Theorem \ref{th4.10}. A straightforward
computation then yields $\pnge\subseteq\cavs$.  Observe that for
$a\in\cavs$, the derivatives $\eta_{r}\partial_{\eta_{s}}\sva$ and
$\partial_{y_{s}}\sva$ of the local symbol of $a$ can be realized as
the local symbol of the commutator of $a$ with one of the operators
mentioned above, hence they satisfy an $L^{2}$ estimate similar to
(\ref{l2est}).  As in \cite{come78} or \cite{st}, an application of
Sobolev's embedding theorem, leads to the following estimate.

\begin{lemma}\label{lemmaest}\
Let $\alpha,\beta\in\NN^{n}$ be arbitrary. Then there exists a
continuous semi-norm $q_{\alpha,\beta}$ on $\cavs$ such that
$$
        \sup_{x,y,\eta}\left|<\eta>^{|\beta|} \partial_{\eta}^{\beta}
        \partial_{y}^{\alpha}\sva(x,y,\eta)\right|\leq q_{\alpha,\beta}(a)\,.
$$ 
\end{lemma}

Let $\cav$ be the closure of $\pnge$ in the submultiplicative
$\Psi^{*}$-algebra $\cavs\cap\ldgch$. By \cite[Corollary 2.5]{aa},
$\cav$ is a $\Psi^{*}$-algebra in $\ldgch$, and $\pnge$ is dense in
$\cav$.  Note that the local symbols of elements in $\cav$ are in the
variable $x$ locally uniform limits of symbols in
$\mathcal{C}(U,S^{0}(T;\rz^{k}_{\eta}))$, hence they remain continuous
in $x$ which gives the following proposition.

\begin{proposition}\label{proplocal}\
Let $a\in\cav$ be arbitrary and 
$\chi:\Omega\stackrel{\cong}{\rightarrow} U\times T$ be as above.
Then for all $\varphi,\psi\in\cunc(\Omega)$ there
exists a symbol $\sigma_{\Omega}(a;\varphi,\psi)\in\mathcal{C}(U,
S^{0}(T;\rz^{k}_{\eta}))$ such that we have
\begin{equation}
\label{locrep}
\chi_{*}(\psi a \varphi)=\sigma_{\Omega}(a;\varphi,\psi)(x,y,D_{y})\,.
\end{equation}
\end{proposition}

Choosing a (countable) cover $\GR=\bigcup_{\Omega\in\cv}\Omega$ of
$\GR$ by open charts $\Omega$ as above, we can define $\can$ to be the
closure of $\pnge$ in the submultiplicative $\Psi^{*}$-algebra
$\bigcap_{\Omega\in\cv}\cav$. Then $\pnge$ is dense in $\can$, and
each $a\in\can$ has a representation (\ref{locrep}) with respect to a
symbol $\sigma_{\Omega}(a;\varphi,\psi)\in\mathcal{C}(U,
S^{0}(T;\rne))$.

As explained above, in a next step we consider now the semi-ideal 
$\cjtu$.

\subsection{Properties of the semi-ideal $\cjtu$} 
The question which properties can be obtained by choosing the set
$\ft$ of closed, symmetric operators appropriately is more complicated
because it includes in particular the analysis of pseudodifferential
operators on non-compact manifolds, hence, we will be rather short at
this point, and sketch only what is within reach.
 
Choose an at most countable set of sections $S\in\cun(M,A(\GR))$,
$S\in\fs$, such that for each $x\in M$ 
the set $\{S(x)\in T_{x}\GR_{x}:S\in\fs_{x}\}$ generates
$T_{x}\GR_{x}$ as a real vector space for some finite subset
$\fs_{x}\subseteq\fs$, and let $\ft\subseteq\pege$ be the corresponding
set of right-invariant, $d$-vertical vector fields. After multiplying them 
with $i$, the elements in $\ft$ have symmetric, minimal closed
extensions by Lemma \ref{lemmamin}. For simplicity, let us assume that 
$\ft$ is finite; otherwise, we have to consider the projective limit
of the corresponding semi-ideals with respect to an increasing
sequence $\ft_{j}\subseteq\ft_{j+1}\subseteq\ft$.
Let $\hft{m}$, $m\in\NN_{0}$ be the corresponding scale of $\ft$-Sobolev
spaces.
Note that the spaces $\hft{m}$ have a decomposition as a direct
integral of the form
$$
        \hft{m}=\dirint_{M}H^{m}_{\ft_{x}}(\GR_{x},\rsed)d\nu(x), 
$$
where $H^{m}_{\ft_{x}}(\GR_{x},\rsed)$ is the Sobolev spaces of order
$m\in\NN_{0}$ associated to the vector fields
$\ft_{x}:=\ft|_{\GR_{x}}\subseteq\pege$. We extend these scales of
Sobolev spaces by duality to $m\in\ZZ$, then the following result
follows immediately from the definition of the ideal $\cjtu_{*}$.

\begin{lemma}\label{lem6}\
Let $a\in\cjtu_{*}$ be arbitrary. Then $a$ induces for each $m\in\NN$
and each $x\in M$ bounded operators $a:\hft{-m}\rightarrow\hft{m}$  and 
$$
        \wh{a}(x):H^{-m}_{\ft_{x}}(\GR_{x},\rsed)\longrightarrow
        H^{m}_{\ft_{x}}(\GR_{x},\rsed).
$$
\end{lemma}

Consequently, $\wh{a}(x)$ has a smooth kernel provided Sobolev's embedding
\begin{equation}\label{eqsobem}
        H^{m}_{\ft_{x}}(\GR_{x},\rsed)\hookrightarrow\cnb(\GR_{x},\rsed)
\end{equation} 
holds for $2m>\dim\GR_{x}$. 
Since everything takes place within the algebra $\can$, the smooth
kernel depends in addition continuously on the parameter $x$. 

Note that (\ref{eqsobem}) holds for instance if the manifolds
$\GR_{x}$ are of bounded geometry \cite[Appendix]{shubin92}.  For the
question, when the groupoid is of bounded geometry we refer to
\cite{aln}.  Let us denote by $\cunb( \GR_{x}\times\GR_{x},
\END(\VD^{\frac{1}{2}}))$ the space of all smooth sections $\GR_{x}
\times\GR_{x} \rightarrow \END(\VD^{\frac{1}{2}})$ that are uniformly
bounded, as are all their covariant derivatives.  A combination of
Lemma \ref{lem6} and \eqref{eqsobem} then yields a characterization of
the Schwartz kernels of the operators $\wh{a}(x)$.

\begin{proposition}\label{probtest}\
Suppose that the manifolds $\GR_{x}$ are of bounded geometry for all
$x\in M$, and let $a \in \cjtu_{*}$ be arbitrary. Then
$\wh{a}(x)\in\lckx$ is an operator with $\cun$-kernel
$k_{\wh{a}}(x,\cdot,\cdot)\in \cunb( \GR_{x} \times \GR_{x},
\END(\VD^{\frac{1}{2}}))$.
\end{proposition}

We summarize the results of the above discussion in the following
proposition.

\begin{proposition}\label{prop.smooth}\
Suppose that for any $x \in M$, the fiber $\GR_x$ is a manifold of
bounded geometry.  Then there exists a $\Psi^*$-algebra $\cae$
containing $\pnge$ as a dense subalgebra such that each $P\in\cae$ is
given by a 
$\GR$-invariant family $(P_x)_{x\in M}$ of
pseudodifferential operators $P_{x}$ on $\GR_x$.
\end{proposition}

Note that the definition of the $\Psi^{*}$-algebra $\cae$ depends on
many choices. 

\begin{proof}\
This follows directly from the plan for the construction of
$\cae$ after Lemma \ref{lemmamin}, Proposition \ref{proplocal} and
Proposition \ref{prop.smooth}. 
\end{proof}

The above results provide us with $\Psi^*$-algebras that are useful in
practice, because they consist of pseudodifferential operators. These
algebras will necessarily contain operators that are not properly
supported (unless our manifold is compact without corners).
Nevertheless, these algebras consist of bounded operators, so their
Schwartz kernels must satisfy some decay conditions far from the
diagonal. It is difficult in general to quantize these decay
conditions.  One possibility is to consider commutators with functions
that approximate the distance function as in
\cite{KordyukovTh}
or  \cite{shubin92}. If the groupoid allows a length function
with polynomial growth as in Section \ref{Sec.6}, we can improve this
by introducing Schwartz spaces.
Also, for certain explicitly given groupoids much more is 
possible; we have elaborated this is Section \ref{Sec.5}
for the case of generalized cusp-calculi on compact manifolds with
corners.

\section{Algebras on manifolds with corners\label{Sec.4}}

In this section we recall the constructions of various groupoids associated to
manifolds with corners (see for instance \cite{lmn,mont97,nwx}). 
We shall use these results to define length functions on some of 
these groupoids, 
which in turn is useful when defining Schwartz spaces associated 
to manifolds with 
corners, in Section \ref{Sec.6}. Note that we do not require that 
the manifolds used
in this section have embedded hyperfaces. Also, we shall use 
these constructions in the
particular case of manifolds with boundary in the next section in order to
construct algebras with smooth kernels. The reader interested only
in the next section, can skip this section at a first lecture, and only refer
back to it when necessary.

Let $X$ be a \mc, and $x$  a point of $X$; we denote by $F(x)$ the connected
component of the set of points  having the same codimension as $x$ which
contains $x$, and by $N_xF(x)=T_xX/T_xF(x)$ the normal space to the boundary 
at $x$.
One can define several groupoids associated to $X$, giving various
pseudodifferential calculi, such as the $b$-calculus, the cusp-calculus, and
its generalizations (the $c_n$-calculi).

Let
$$
        \GG(X)\!\!=\!\!\Big\{ (x,y,\alpha)\,|\, x,y\! \in\! X,\; 
        {\rm codim}(x)\!\!=\!\!{\rm codim}(y),\;
        \alpha\!:\!N_yF(y) \build{\to}_{}^{\sim} N_xF(x) \Big\},
$$
where $\alpha$ is given, through  
trivializations  $N_yF(y) \simeq \r+{k}$ based on inward pointing normals
and $N_xF(x) \simeq \r+{k}$, by a matrix which has one and only one non-zero
element on each line and each column, and this element is positive. It is
precisely the product of a diagonal matrix with all  terms strictly
positive by a permutation matrix.  

The groupoid structure of $\GG(X)$ is given by $d,r : \GG(X) \to X$ with
$r(x,y,\alpha)=x$, and $d(x,y,\alpha)=y$; the composition law is induced by the
composition of the isomorphisms, in the sense that $(x,y,\alpha)(y,z,\beta)
=(x,z,\alpha \beta )$.

One can endow $\GG(X)$ with several different differential structures
such that the resulting groupoids are homeomorphic but not
diffeomorphic.  Let $\Omega$ and $\Omega'$ be two charts of $X$ of
same codimension. Thus $\Omega \build{\simeq}_{}^{\pi} U \times
\r+{A}$ and $\Omega'\build{\simeq}_{}^{\pi'} U' \times \r+{A'}$, where
$U$ and $U'$ are open subsets of some $\R^n$ and $A$ and $A'$ are the
sets of local hyperfaces contained in $\Omega$ and $\Omega'$. The sets
$A$ and $A'$ have the same cardinal, namely the codimension of
$\Omega$. Fix a bijection $\sigma : A' \to A$.

If $t \in \r+{A'}$ and $\lambda \in \rep{A'}$, then let $B'_t=\{ i \in
A', t_i=0\}$. The product of the matrix of $\sigma_{|B'_t}$ by the
diagonal matrix consisting of the $\lambda_i$, for $i\in B'_t$, is
denoted by $m_{\sigma,t,\lambda}$ ; it defines an isomorphism
$\alpha_{\sigma,t,\lambda}$ through the trivializations induced by
$\pi$ and $\pi'$:
$$
        \begin{array}{ccc}
        \alpha_{\sigma,t,\lambda} : N_{{\pi'}^{-1}(u',t)}F({\pi'}^{-1}(u',t))
        &\hfl{\simeq}{} &N_{\pi^{-1}(u,\sigma(\lambda
        t))}F(\pi^{-1}(u,\sigma(\lambda t)))\\
        \vfl{\simeq}{} & & \vfl{\simeq}{} \\
        \r+{B'_t} & \hfl{m_{\sigma,t,\lambda}}{} & \r+{B}
        \end{array}
$$

\subsection{The $b$-calculus differential structure}
Using the notations above, consider the map
$$
        \begin{array}{rcl}
        \psi_\sigma : U \times U' \times \r+{A'} \times \rep{A'} &\to& \GG(X) \\
        (u,u',t,\lambda) & \mapsto & (\pi^{-1}(u,\sigma(\lambda t)), {\pi'}^{-1}(u',t), 
        \alpha_{\sigma,t,\lambda}).
        \end{array}
$$

This map is injective, and one thus obtains a \cinf\ structure on $\GG(X)$,
which turns it into a Lie groupoid, whose fibers are \sm s without
boundary; it is amenable.

\begin{definition}\
  The groupoid of the \bcalc, 
$\Gamma_1(X)$ is the union of the connected
components containing the unit of each $d$-fiber of $\GG(X)$.
\end{definition}

\begin{remark}\
It is instructive to consider the special case when $X$ is a manifold with 
connected boundary $\pa X$ and defining function $\rho$. 
Then one can prove that 
$$
        \Gamma_1(X)\simeq\{ (x,y,\lambda) \in X\times X\times \rep{}\,|\, 
        \rho(x)=\lambda \rho(y) \}.
$$
The identification with the $b$-calculus is obtained by
observing that $\pa M \times \pa M \times \rep{} \subset \Gamma_1(X)$ and
that $\log \lambda = \log \rho(x) - \log \rho(y)$ if $(x,y,\lambda) \in
\Gamma_1(X)$ and $x \not \in \pa M$.
\end{remark}

\subsection{The cusp-calculus and $c_n$-calculi differential structures}
Consider a continuous, strictly increasing map,
smooth for $t \not = 0$, such that
\begin{equation*}
        \begin{split}
        \tau_n :& \r+{} \to \r+{} \\
        & t \mapsto 
        \begin{cases}
      \frac{1}{e}(-\log(t))^{-\frac{1}{n}} & \mathrm{ if\ } t\in (0,1/e)\\
      0 & \mathrm{ if\ } t= 0\\
      t & \mathrm{ if\ } t\geq 1
      \\
        \end{cases}
        \end{split}
\end{equation*}

If $U$ is an open subset of a Euclidean space, we will also 
denote by $\tau_n : U \times\r+{A'} 
\to U\times\r+{A'}$ the map obtained by
applying $\tau_n$ to each coordinate of $\r+{A'}$. 

As above, one can define
$$
        \begin{array}{rcl}
        \psi_{\sigma,n} : U \times U' \times \r+{A'} \times \rep{A'} &\to& 
           \GG(X) \\
        (u,u',t,\lambda) & \mapsto & (\pi^{-1}(u,\sigma(\tau_n(\lambda t))), 
        {\pi'}^{-1}(u',\tau_n(t)), \alpha_{\sigma,t,\lambda})
        \end{array}
$$
which endows $\GG(X)$ with a new differential structure.

\begin{definition}\ Let $n \ge 2$.
The groupoid of the $c_{n}$-calculus, $\Gamma_{n}(X)$, 
is the union of the connected
components containing the unit of each $d$-fiber of the \gr\ $\GG(X)$
endowed with the structure defined by 
$\psi_{\sigma,n-1}$ (using $\tau_{n-1}$). 
\end{definition}

Recall that the $c_2$-calculus is also known as the \textit{cusp-calculus}
-- see for instance \cite{in,MazzeoMelrose,MelroseFB,meni96c}

\begin{remark}\
When $X$ is a manifold with connected
boundary, endowed with a defining function of the boundary, $\rho$, then
$$
        \Gamma_{n+1}(X)\simeq\{ (u,v,\mu) \in X\times X\times \R\,|\, 
        \mu \rho(u)^n\rho(v)^n = \rho(u)^n-\rho(v)^n\}
$$
as smooth manifolds, which can be seen directly from the definition.
Moreover, the structural morphisms of the groupoid $\Gamma_{n+1}$ become
$d(u,v,\mu) = v$, $r(u,v,\mu) = u$, and $(u,v,\mu) (v,w,\lambda) =
(u,w,\mu + \lambda)$. The Lie algebroid of this groupoid is seen to
consist of the vector fields $X \in \Gamma(TM)$ such that $X(\rho(x)^{-n})$
is a smooth function on $M$.
\end{remark}

\subsection{Comparison of the $c_n$-pseudodifferential calculi}

The groupoids defined above only differ by their differential
structures; in fact, they are all homeomorphic. This is intuitively
clear since $\tau_n$ is a homeomorphism (but it is NOT a
diffeomorphism) which induces a homeomorphism of groupoids. An
immediate application is the fact that the $b$-calculus and
$c_n$-calculi have the same norm closure.

To keep notations simple, we will only consider here the case of a
manifold with boundary.  Then by considering a collar neighborhood of
the boundary, one gets a partition
$$
        X=X_1 \cup X_2
$$
with $\pi : X_1 \simeq \dx \times [0,1)$ and 
$X_2=X\setminus \pi^{-1}( \dx \times [0,1/e))$; 
the boundary defining function used here is 
$$    
        \rho(x)=        \begin{cases}
        e  t & \mathrm{\ if \ } x \in X\setminus X_2 \mathrm{\ with\ } 
        \pi(x)=(u,t)\\
        1 & \mathrm{\ if \ } x \in X_2\\
        \end{cases}
$$
($\rho$ is not a smooth function on $X$, but it is smooth on $\dx$).

This allows us to define a homeomorphism

\begin{equation}\label{theta}
  \begin{split}
  \Theta_{n+1} :  \Gamma(X) &\to \Gamma_{n+1}(X)\\
(x,y,\lambda) & \mapsto (u,v,\mu)
  \end{split}
\end{equation}
with 
\begin{equation*}
  \begin{split}
    u=&
\begin{cases}
  x \mathrm{\ if\ } x \in X_2\\
  \pi^{-1}\circ \tau_n \circ \pi(x) \mathrm{\ if\ } x \in X\setminus X_2
\end{cases}\\
v=&
\begin{cases}
  y \mathrm{\ if\ } y \in X_2\\
  \pi^{-1}\circ \tau_n \circ \pi(y) \mathrm{\ if\ } y\in X\setminus X_2
\end{cases}\\
\mu=&\log(\lambda)
\end{split}
\end{equation*}

One can check then that $\mu \rho(u)^n\rho(v)^n=\rho(u)^n-\rho(v)^n$ 
by considering
separately the following cases:
\begin{enumerate}
\item If $x \in X\setminus X_2$ and $y \in X\setminus X_2$ then $u \in
        X\setminus X_2$ and $v \in X\setminus X_2$, since $\tau_n(t)\leq 1/e$ if
        $t\leq 1/e$. Thus if $\rho(u)\neq 0$, $\rho(u)^{-n}=-\log(\lambda t)$ 
         and
        $\rho(v)^{-n}=-\log(t)$, which implies that
        $\rho(v)^{-n}-\rho(u)^{-n}=\log(\lambda)$. If $\rho(u)=0$, 
        then $\rho(v)=0$
        and the equality is trivial.
\item If $x \in X\setminus X_2$ and $y \in X_2$, then $u \in
        X\setminus X_2$ and $v=y \in X_2$. Thus 
        $\rho(u)^{-n}=-\log(\rho(x)/e)$, and $\rho(v)^{-n}=1$. But
        $\rho(x)=\lambda\rho(y)=\lambda$, so that
        $\rho(v)^{-n}-\rho(u)^{-n}=\log(\lambda)$. 
\item If $x \in X_2$ and $y \in X_2$, then $u=x \in
        X_2$ and $v=y \in X_2$. Thus $\rho(u)=\rho(v)=1$ and the equality is
        trivial. 
\end{enumerate}

\section{Algebras of smooth kernels\label{Sec.5}}

In this section we shall study regularizing operators on certain 
differentiable groupoids on a manifold with boundary.

Recall that a {\em differentiable groupoid} is a continuous family groupoid 
$\GR$ such that the space of units $M$, as well as the space of arrows
$\GR^{(1)}=\GR$, are differentiable manifolds (possibly with corners),
all structural map are differentiable, and the domain map
$d : \GR\rightarrow M$ is a submersion of manifolds with corners.  Note
that the latter in particular implies that the fibers
$\GR_{x}:=d^{-1}(x)$ are smooth manifolds without corners that are in
general non-compact.  Pseudodifferential operators on differentiable
groupoids have been considered in \cite{sur,mopi97}, and \cite{nwx} in
more detail.

We shall use now the results of the previous sections to construct an
algebra of regularizing operators that is closed under holomorphic
functional calculus and whose kernels are smooth {\em including on the
boundary}, for suitable $\GR$. This is non-trivial, in view of the
results of \cite{aa,Melrose42}, where it is proved that this is not
possible for the $b$-calculus. We begin by formulating the problem
more precisely.

Let $\GR \to M$ be a Hausdorff differentiable groupoid on a manifold
with corners $M$. We want to construct algebras $A$ with the following
properties:

\begin{enumerate}
\item $\Psi^{-\infty}(\GR) = \CIc(\GR) 
      \subset A \subset \CI(\GR) \cap C^*(\GR)$ and
\item $A$ is a (possibly non-unital) $\Psi^*$-algebra.
\end{enumerate}

\begin{definition}\ An algebra $A$ satisfying properties (a) and (b)
right above is called {\em  a $\Psi^*$-algebra of smooth kernels on $\GR$}. 
\end{definition}

In \cite{aa} it is
proved that there is no $\Psi^*$-algebra of smooth kernels on
$\GR=\Gamma_1(M)$, where $\Gamma_1(M)$ is the smooth groupoid
associated to the $b$-calculus (see Section \ref{Sec.4}). However, we
shall now show how to construct algebras of smooth kernels on $\GR$,
if $\GR = \Gamma_n(M)$ are the groupoids defining the $c_n$-calculi on
a manifold with boundary $M$, provided that $n \ge 2$ (see Section
\ref{Sec.3} for the definition of $\Gamma_n(M)$). 
Some of us have  learned that it is
possible to construct algebras of smooth kernels on $\Gamma_n(M)$
from Richard Melrose.

Let $I := \dot{\mathcal{C}}^{\infty} (M \times M)$ 
be the space of smooth functions on $M \times M$ that vanish to
infinite order on the boundary (the boundary here is the union of
hyperfaces of $M \times M$).  Then $\dot{\mathcal{C}}^{\infty} (M
\times M)$ is an algebra of smooth kernels, and hence an algebra of
order $-\infty$ pseudodifferential operators on $M_0 := M
\smallsetminus \pa M$.

Moreover, these operators are bounded on $L^2(M_0)$. We have to note
here that the correct density on $M_0$ for the $c_n$-calculus is, on a
tubular neighborhood of $\pa M$, of the form $hx^{-n}|dx||dy|$, where
$x \in [0,\infty)$ is a boundary defining function on $M$ and $|dy|$
is a density on $\pa M$, and $h \in \CI(M)$.

\begin{lemma}\label{lemma.I}\ The space $I \subset C^*(\GR)$ is a 
non-unital $\Psi^*$-algebra.
\end{lemma}

\begin{proof}\ Let $\Delta$ be
the Laplace operator for some compatible metric (i.e. such that
$x^n\pa_x$ has length one). 
Take $\ft = \{\Delta, x^{-1}\}$, 
regarded as unbounded operators on $L^2(M)$, and apply the semi-ideal
construction to $C^*(\Gamma_n(M))$ and $\ft$. This yields an algebra
$J:=\cjt{-\infty}_{*}$ satisfying
$$
        J \subset \{ T \in C^*(\Gamma_n(M))\, |\;
        x^{-i}\Delta^{j}T\Delta^{k}x^{-l}\, \text{ is bounded }\,
        \forall i,j,k,l\,\}.
$$
Clearly, the opposite inclusion is also true by the definition of
$\cjt{-\infty}$.

Let $H^m(M)$ be the domain of $\Delta^{m/2}$, if $m \ge 0$, or 
$H^m(M)$ be the dual of $H^{-m}(M)$, if $m$ is negative. Also, let
$H^{-\infty}(M) = \cup H^m(M)$ and $H^\infty(M) = \cap H^m(M)$.

Since every $T \in J$ maps $H^{-\infty}(M) \to H^{\infty}(M)$ and $J$
is symmetric, we obtain that $T$ is a smoothing operator. Thus, $T$ is
an operator with integral kernel given by a smooth function $K(x,y)$.
Since
$x^{-i}K(x,y,x',y'){x'}^{-j}$ also must define a bounded operator, we
see that $K$ vanishes to infinite order at the boundary.
Consequently, $J=I$, and Corollary \ref{cor.psi.J} completes the
proof.
\end{proof}

\begin{remark}\label{rem.alt}\ 
In the  proof above, we  obtain the same conclusion 
by considering $\ft = \{ x^{-1}, X_1, \ldots, X_m\}$, where 
$X_1, \ldots, X_m \in \Gamma(TM)$ is a system of generators for the 
vector fields corresponding to the $c_n$-calculus, that is, 
$b$-vector fields satisfying $X_j (x^{-n + 1}) \in \CI(M)$.
\end{remark}

We proceed now to describe the regularizing operators in the
$c_n$-calculi on a manifold with boundary $M$ \cite{sur} in a way that
is most convenient for our purposes.

Let $(\ca,(\norm{\cdot}_{j})_{j\in\NN})$ be a submultiplicative
Fr\'echet algebra. Assume there is given an action $\alpha : \RR \to
\operatorname{Aut}(\ca)$ of $\RR$ by automorphisms on $\ca$ (so
$\alpha_t \circ \alpha_s = \alpha_{t+s}$). If, for any $a \in \ca$,
the map $\RR \ni t \mapsto \alpha_t(a)$ is smooth and there exist
polynomials $P_{n}$, 
{\em independent of $n$}, 
such that
\begin{equation}
        \norm{\alpha_t(a)}_n \le P_{n}(|t|) \norm{a}_n,
\end{equation}
then we say that the action of $\RR$ on $\ca$ is {\em with polynomial
growth}. In that case, we can introduce on the Schwartz space
$\SS(\RR, \ca)$ an algebra structure by
$$
        f * g (t) = \int_\RR f(s) \alpha_s(g(t - s)) ds.
$$
Moreover, $\SS(\RR,\ca)$ acquires a family of seminorms $\|\;\|_{n,i,j}$,
$$
        \|f\|_{n,i,j} = \int_{t \in
        \RR}\|t^i\pa_t^jf(t)\|_ndt,
$$
making it a Fr\'echet algebra, that is submultiplicative with
respect to an {\em equivalent} system of
seminorms that we now define.

First, we can assume that 
$P_{n}(t) = C_{n}(1 + t^{M_{n}})$, for
some $C_{n} > 0$ and $M_{n} \in \NN$. Then, a direct verification
using the binomial expansion of $s^i = (s - t + t)^i$ and the
submultiplicativity of the seminorm $\| \; \|_n$ of $\ca$ shows that
\begin{multline}
\label{eq.robert}
        \|f * g \|_{n,i,j} = \int \| t^i \pa_t^j\int f(s)
        \alpha_{s}(g(t-s)) ds \|_n dt \\
        \le C_{n}\sum_{\beta + \gamma = i} C_i^\beta \big(
        \|f\|_{n,\beta,0} + \|f\|_{n,\beta+M_{n},0} \big) \| g
        \|_{n,\gamma,j}.
\end{multline}
(Here $C_a^b = a!b!^{-1}(a-b)!^{-1}$ are the binomial coefficients.)

Let us adjoin a unit denoted $e$ to $\SS(\RR,\ca)$ and extend the norms $ \|\,\cdot\,\|_{n,i,j} $
to $\SS(\RR,\ca)^+ := \SS(\RR,\ca) \oplus \CC e$ by
$\| (f,\lambda e) \|_{n,i,j} := \| f \|_{n,i,j} + |\lambda|$.
Equation \eqref{eq.robert} implies that left multiplication by $f$
on $\SS(\RR,\ca)^+$ is continuous in the norm
$\sum_{\alpha \le i} \| \, \cdot \, \|_{n,\alpha,j}$. The corresponding
operator norm, $|||f|||_{n,i,j}$ is then submultiplicative, by
definition and satisfies
$$
        |||f|||_{n,i,j} \le \|f\|_{n,i,j} + K_{n,i,j} \sum_{l \le i + M_n} \|f\|_{n,l,0} \,.
$$
(The first term appears due to the fact that we have adjoined a unit to $\SS(\RR,\ca)$.
Also,  $K_{n,i,j}$ is a constant that is independant of $f$.)
On the other hand, 
$$
        \|f\|_{n,i,j} = \| f * e \|_{n,i,j} \le |||f|||_{n,i,j} \|e\|_{n,i,j}
        = |||f|||_{n,i,j}, 
$$
which proves that the families of seminorms $|||\,\cdot\,|||_{n,i,j}$ and 
$\|\,\cdot\,\|_{n,i,j}$ define the same topology on $\SS(\RR,\ca)$.

Let $M$ be a compact manifold with boundary $\pa M$.  On $\pa M \times
[0,\infty)$ we consider the vector field $X_n = (1 + x^n)^{-1} x^n \pa
_x$, with $x \in [0,\infty)$. Our choice of this vector field is
justified by the fact that 
$$
        X_n(x^{-n+1}) = (-n + 1) (1 + x^n)^{-1} \in \CI([0,\infty)),
$$
and hence it defines a vector field corresponding to the $c_n$-calculus, 
 $n \ge 2$. Moreover, this vector field defines, by
integration, an action of $\RR$ on $\pa M \times \pa M \times
[0,\infty)$, which is trivial on $\pa M$. 

Let $\ca := \SS(\pa M \times \pa M \times [0,\infty))$, with product 
$$
        (f g)(m_1,m_2,t) = \int_{\pa
        M} f(m_1,m,t) g(m,m_2,t) dm.
$$
Then $\ca$ is isomorphic, as an algebra, with the 
complete projective tensor product
$\Psi^{-\infty}(\pa M) \otimes_\pi \SS([0,\infty))$.

\begin{lemma}\label{lemma.pol.gr}\
Assume $n \ge 2$. Then the action of $\RR$ on $\ca := 
\SS(\pa M \times \pa M \times [0,\infty))$ integrating the action of the 
vector field $X_n = (1 + x^n)^{-1} x^n \pa _x$ is with polynomial growth.
\end{lemma}

\begin{proof}\
Let $S^0(\RR)$ be the space of classical symbols of order zero
on $\RR$, with its natural Fr\'echet space structure. Consider the function 
\begin{equation}
        f_n(x) = (-n + 1)^{-1} x^{-n + 1} + x : (0,\infty) \to \RR.
\end{equation}
Then $f_n$ is a bijection such that the induced map 
$$
        f_n^*:\SS([0,\infty)) \to \sum_{k = 0}^{n-2}S^{-k/(n-1)}(\RR)
$$ 
is an equivariant map, i.e.\ the action of $\RR$ on $\SS(\RR)$ being by
translation. Moreover $f_n^*$ is an isomorphism onto its image, which
consists of the linear combinations of symbols that are of order
$-\infty$ on the positive semi-axis $[0,\infty) \subset \RR$.  In
particular, the image of $f_n^*$ is closed.

Since the action of $\RR$ by translation on $S^k(\RR)$ is
with polynomial growth, the given action of $\RR$ on $\SS([0,\infty))$ 
is also with polynomial growth.
\end{proof}

Probably the above lemma is the only reason why we have to restrict to
$n \ge 2$ in order to construct algebras with smooth kernels on
the $c_n$-calculus groupoid $\Gamma_n(M)$.

Fix a smooth function $\phi \in \CI([0,\infty))$, $\phi(x) = 1$ if $x$
is in a certain neighborhood of $0$, $\phi(x) = 0$ if $x \ge 1$, and
let
$$
        A_1 = \phi \SS(\RR, \ca) \phi := \phi \SS(\RR, \SS(\pa M
        \times \pa M \times [0,\infty)))\phi,
$$
where $\ca := \SS(\pa M \times \pa M \times [0,\infty))$, as in the
above lemma. The algebra $A_1$ acts on $L^2(\pa M \times
[0,\infty))$. 

Let $\Gamma_n(M)$ be the groupoid defining the $c_n$-calculus. Then
$\Psi^{-\infty}(\Gamma_n(M))$ identifies with an algebra
of bounded operators on $L^2(M)$ (with the canonical induced
measure).

We are ready to prove the following theorem.

\begin{theorem}\label{theorem.smooth.b}\ Let $M$ be a compact, smooth
manifold with connected boundary.
Then $\mathfrak A = A_1 + I \subset C^*(\Gamma_n(M))$
and is a non-unital $\Psi^*$-algebra with smooth kernels.
\end{theorem}

\begin{proof}\ We shall use the results of the previous sections. First, 
however, we have to prove that $\mathfrak A := A_1 + I $
is an algebra.

Indeed, by Theorem~\ref{thm4.12}, it is enough to check that 
$A_1 \subset \puft$, 
where $\ft = \{x^{-1},X_0,\ldots,X_m\}$ is 
chosen as in Remark \ref{rem.alt}.
We identify a tubular neighborhood of $\pa M$ with a subset of
$\pa M \times [0,\infty)$. To prove this, we first notice that
$X_0,\ldots,X_m \in \Psi^{1}(\Gamma_n(M))$, 
$A_1 \subset \Psi^{-\infty}(\Gamma_n(M))$,

and hence any  product of the form 
$$
        X_{i_1} \ldots X_{i_j} a X_{i_{j+1}} \ldots X_{i_r}
$$
consists of bounded operators.  Then, let us write by $a_f$ the
operator of convolution on $\pa M \times [0,\infty)$ with the function
$f \in \SS(\RR)$ (the action of $\RR$ is the one integrating the
vector field $X_n = (1 + x^n)^{-1}x^n \pa_x$). Also, let us observe
that $\ca = \SS (\pa M \times \pa M \times [0,\infty))$ identifies
with an algebra of operators on $L^2(\pa M \times [0,\infty))$.  Then
the map
$$
        \SS(\RR) \otimes_\pi \ca \ni f \otimes b \to a_f b \in
        \SS(\RR,\ca)
$$
is an isomorphism. We need to check that $\delta$, $\delta(a) :=
[x^{-1},a]$ maps $\SS(\RR,\ca)$ to itself. Clearly $\delta(b) = 0$,
for any $b \in \ca$. If $f \in \SS(\RR)$, then $\delta(f) \in
\SS(\RR,\ca)$.  This proves our claim that $\SS(\RR,\ca)$ is stable
with respect to $\delta$. In conclusion, $A_1 \in \puft$ and
hence $A_1I + IA_1 \subset I$ (see Theorem \ref{thm4.12}).

The algebra $\mathfrak A$ has a Fr\'echet topology induced from the
Frechet topologies of $A_1$ and $I$. To prove that it is a
$\Psi^*$-algebra, we thus only need to prove that it is spectrally
invariant. To this end, we shall use the results of Theorem
\ref{thm4.02}. Let $\cb = C^*(\Gamma_n(M))$, regarded as a subalgebra
of the algebra of bounded operators on $L^2(M)$. Also, let $\cj$ be
the norm closure of $x\cb$. Let $\mathcal K$ be the algebra of compact
operators on $L^2(\pa M)$.  Then $\cb/\cj \simeq C_0(\RR,\mathcal K)$
by standard results on groupoid $C^*$-algebras. (This statement can
also be extracted from either of \cite{lmn,MelroseNistor}, by using,
for example the exact sequence associated to the invariant subset $\pa
M \times \{0\} \subset \pa M \times [0,\infty)$.) Let $\varphi :
\mathfrak A/ I \to \cb/\cj$ be the induced morphism. Then the range of
this morphism is $\SS(\RR, \CI(\pa M \times \pa M)) \subset \cb/\cj
\simeq C_0(\RR, \mathcal K)$
and hence $\varphi$ is locally spectral invariant.

Theorem \ref{thm4.02} or, more precisely Corollary \ref{coreins}, then
shows that $\mathfrak A$ is spectrally invariant. This completes the
proof.
\end{proof}

We are planning to clarify the relation between our algebra
$\mathfrak{A}$ for the $c_{2}$-calculus and the cusp-calculus as
defined by Melrose (an exposition can be found in \cite{meni96c}) in
another paper.

\section{The Schwartz space of a continuous family groupoid\label{Sec.6}}

In this section we define a notion of Schwartz space on a \cfg\ $\GR$,
i.e.\ a space of rapidly decreasing functions as well as their
derivatives. This was introduced in \cite{mont98} in the case of
differentiable groupoids. We prove, using also some methods introduced
in \cite{lafforgue-these}, that this is a subalgebra of $C^*(\GR)$,
stable under holomorphic functional calculus.

\subsection{The Schwartz convolution algebra}
Let $\GR$ be a Hausdorff, \cfg. Fix a $1$-density on $A$. This will
then give rise to a $1$-density on each of the manifolds $\GR_x$, and
hence to a smooth measure $\mu_x$ on each of $\GR_x$. Moreover, the
measures $\mu_x$ are invariant with respect to right translations, and
hence they form a Haar system.

\begin{definition}\
Let $\mu$ be the Haar system on $\GR$ introduced above.  A length
function with polynomial growth on $\GR$ is a continuous function
$\phi~: \GR \to \R_+$ such that:
\begin{enumerate}
\item $\phi(g_1g_2)\leq 
\phi(g_1)+\phi(g_2)$, 
\item $\forall g \in \GR, \phi(g^{-1})=\phi(g)$,
\item $\phi$ is proper,
\item  $\exists c, N, \forall x \in \GR^{(0)}, \forall r \in \R_+, \mu_x(\phi^{-1}([0,r]) \leq c (r^N+1)$.  
\end{enumerate}
\end{definition}

The polynomial growth condition ensures that there exists $k_0 \in
\NN$ and $C \in \RR$ such that for any $k \geq k_0$, and for any $x\in
\g0$,
$$ 
        \int_{G_x} \frac{1}{(1+\phi(g))^{k}}d\mu_x \le C.
$$

Let $v$ be a section of the Lie algebroid of $\GR$, $A(\GR)$; such a
section defines a differential operator of order 1 on $\GR$. Thus if
$v_1, \ldots, v_l$ are sections of $A(\GR)$, and if $f \in {\cal
C}_0(\GR, \O12)$, then $v_1 \ldots v_k \cdot f \cdot v_{k+1} \ldots
v_l$ is a distribution on $\GR$. It belongs to ${\cal C}_0(\GR, \O12)$
provided that there exists $g \in {\cal C}_0(\GR, \O12)$ such that,
for any $a \in \ccinfo(G)$, $g \cdot a=(v_1 \ldots v_k \cdot f \cdot
v_{k+1} \ldots v_l) \cdot a$.

\begin{definition}\
Let $\GR$ be a \cfg\ and $\phi$ be a length function with polynomial
growth on $\GR$.  Define
\begin{multline*}
        \SS^{k,d}(\GR, \phi)=\big\{ f \in
        \mathcal{C}_0(G,\Omega^{\frac{1}{2}}), \forall v_1, \ldots,
        v_d \in \mathcal{C}(A(\GR)), \forall i \leq d,\\ v_1 \ldots
        v_i \cdot f \cdot v_{i+1} \ldots v_d \in {\cal C}_0(\GR, \O12)
        \\ \mathrm{\ and \ } \sup_{g \in \GR} |v_1 \ldots v_i \cdot f
        \cdot v_{i+1} \ldots v_d(g)|(1+\phi(g))^{k} < \infty \}.
\end{multline*}
The Schwartz space of $\GR$ with respect to $\phi$ is
$$
        \SS(\GR, \phi)=\bigcap_{k,d \in \NN} \SS^{k,d}(\GR, \phi).
$$
\end{definition}

The space $\SS^{k,d}(\GR, \phi)$ can be endowed with the  norm
$$
        \n{f}{k,d}=\sup_{i\leq l\leq d} \sup_{\substack{ v_1, \ldots,
        v_l \in \mathcal{C}(A(\GR))\\ \n{v_j}{}\leq 1}} \sup_{g \in
        \GR} |v_1 \ldots v_i \cdot f \cdot v_{i+1} \ldots
        v_l(g)|(1+\phi(g))^{k}.
$$

\begin{remark}\label{remark}
\begin{itemize}
\item If $l\geq k$ and $f \in \SS^{l,d}(\GR, \phi)$ then $f \in
\SS^{k,d}(\GR, \phi)$ and $\n{f}{k,d}\leq \n{f}{l,d}$.
\item If $f \in \SS^{k,d}(\GR, \phi)$ and $k\geq k_0$ then for any $x\in \g0$ and any $v_1,
\ldots, v_l \in \mathcal{C}(A(\GR))$ (with $l\leq d$), one has
\begin{equation*}
\begin{split}
        \int_{\GR_x} |v_1 \ldots v_i \cdot f \cdot v_{i+1} \ldots
        v_l(g)|^2 &\leq \n{f}{k,d}^2 \int_{\GR_x}
        \frac{1}{(1+\phi(g))^{2k}} d\mu_x\\ &\leq C \n{f}{k,d}^2
\end{split} 
\end{equation*} 
so that $v_1 \ldots v_i \cdot f \cdot v_{i+1} \ldots v_l\in
L^2(\GR_x)$, and $\n{v_1 \ldots v_i \cdot f \cdot v_{i+1} \ldots
v_l}{L^2(\GR_x)}\leq \sqrt{C} \n{f}{k,d}$.
\end{itemize} 
\end{remark}

\begin{proposition}\label{subalgebra}\
For any $k, d \in \NN$, with $k$ such that $\int_{\GR_x}
\frac{1}{(1+\phi(g))^{k}}d\mu_x$ converges for any $x \in \g0$,
$\SS^{k,d}(\GR, \phi)$ is a dense subalgebra of $\cred{\GR}$, and
there exists a constant $\lambda_k$ such that $\n{\cdot}{\cred{\GR}}
\leq \n{\cdot}{k,d}$.
\end{proposition}

\begin{proof}\
If $f_1, f_2 \in \SS^{k,d}(\GR, \phi)$. Let $x = d(g)$, we then have
\begin{equation*}
\begin{split}     
        &|v_1 \ldots v_i \cdot f_1 * f_2\cdot v_{i+1} \ldots v_l(g)| =
        |\int_{\GR_{x}} 
        v_1 \ldots v_i \cdot f_1(g{g'}^{-1})f_2\cdot v_{i+1}  \ldots
        v_l(g')| d\mu_x\\ 
        \leq & \int_{\begin{subarray}{l}
        \GR_{x}\\
        \phi(g')\geq \frac{\phi(g)}{2}
        \end{subarray}}
        | v_1 \ldots v_i \cdot
        f_1(g{g'}^{-1})f_2\cdot v_{i+1}  \ldots v_l(g')| d\mu_x\\
        &+ \int_{\begin{subarray}{l}
        \GR_{x}\\
        \phi(g{g'}^{-1})\geq  \frac{\phi(g)}{2}
        \end{subarray}}
        |v_1 \ldots v_i \cdot
        f_1(g{g'}^{-1})f_2\cdot v_{i+1}  \ldots v_l(g')| d\mu_x\\  
        \leq & \int_{\begin{subarray}{l}
        \GR_{x}\\
        \phi(g')\geq \frac{\phi(g)}{2}
        \end{subarray}}
        \n{f_1}{k,d}\frac{1}{(1+\phi(g{g'}^{-1}))^k} 
        \n{f_2}{k,d}\frac{1}{(1+\phi(g'))^k} d\mu_x\\
        &+  \int_{\begin{subarray}{l}
        \GR_{x}\\
        \phi(g{g'}^{-1})\geq \frac{\phi(g)}{2}
        \end{subarray}} 
        \n{f_1}{k,d}\frac{1}{(1+\phi(g{g'}^{-1}))^k}
        \n{f_2}{k,d}\frac{1}{(1+\phi(g'))^k} d\mu_x\\ 
        \leq & \frac{2^k}{(1+\phi(g))^k} \n{f_1}{k,d}\n{f_2}{k,d}
        \biggl(\int_{\GR_{x}} \frac{1}{(1+\phi(g{g'}^{-1}))^k}+
        \int_{\GR_{x}} \frac{1}{(1+\phi(g'))^k}\biggr) d\mu_x\\
        \leq & \frac{2^{k+1}C}{(1+\phi(g))^k} \n{f_1}{k,d}\n{f_2}{k,d}
        \end{split}
\end{equation*}
so that 
\begin{equation}\label{convolution}
        |v_1 \ldots v_i \cdot f_1 * f_2\cdot v_{i+1} \ldots
        v_l(g)|(1+\phi(g))^k \leq 2^{k+1}C\n{f_1}{k,d}\n{f_2}{k,d}
\end{equation}
which implies that $f_1*f_2 \in \SS^{k,d}(\GR, \phi)$. Consequently,
$\SS^{k,d}(\GR, \phi)$ is an algebra.

To prove that $\SS^{k,d}(\GR, \phi)$ is a subspace of $\cred{\GR}$, we
need to show that if $f\in \SS^{k,d}(\GR, \phi)$, and for any $x\in
\g0$ and $\xi \in C_c(\GR_x)$, one has
$$
        \n{f*\xi}{L^2(\GR_x)} \leq \n{f}{k,d} \n{\xi}{L^2(\GR_x)}
$$
up to a constant.

Denote by $\lambda_k$ the constant $2^{k+1}C$. Then the Cauchy-Schwarz
inequality implies that
$$
        \n{f\xi}{L^2(\GR_x)}^2 \leq
        \n{\xi}{L^2(\GR_x)}\n{(f^*f)\xi}{L^2(\GR_x)}.
$$
By induction, 
$$
        \n{f\xi}{L^2(\GR_x)}^{2^{n+1}} \leq
        \n{\xi}{L^2(\GR_x)}^{2^{n+1}-1}\n{(f^*f)^{2^n}\xi}{L^2(\GR_x)}.
$$
But if $h\in \SS^{k,d}(\GR, \phi)$,
\begin{equation*}
\begin{split}
        \n{h\xi}{L^2(\GR_x)}^2 & = \int_{\GR_x}
        (h(g{g'}^{-1})\xi{g'})^2 d\mu_x\\ & \leq \n{h}{k,d}^2
        \int_{\GR_x} (\frac{\xi{g'}}{(1+\phi(g{g'}^{-1}))^k})^2 d\mu_x
\end{split}
\end{equation*}

The inequality \eqref{convolution} also gives
$$
        \n{(f^*f)^{2^n}}{k,d} \leq \lambda_k^{2^{n+1}} \n{f^{2^{n+1}}}{k,d},
$$ 
so that
$$
        \n{f\xi}{L^2(\GR_x)}^{2^{n+1}} \leq
        \n{\xi}{L^2(\GR_x)}^{2^{n+1}-1} \lambda_k^{2^{n+1}}
        \n{f}{k,d}^{2^{n+1}} (\int_{\GR_x}
        (\frac{\xi{g'}}{(1+\phi(g{g'}^{-1}))^k})^2)^{1/2} d\mu_x.
$$
Taking $n\to \infty$,
$$
        \n{f\xi}{L^2(\GR_x)}\leq \n{\xi}{L^2(\GR_x)} \lambda_k
        \n{f}{k,d},
$$
which implies that $f \in \cred{\GR}$, and $\n{f}{\cred{\GR}} \leq
\lambda_k\n{f}{k,d}$.
\end{proof}

\subsection{The theorem} 
We are now ready to prove the main result of this section, the fact
that the algebras $\SS(\GR,\pi)$ are closed under holomorphic
functional calculus.

\begin{theorem}\label{theorem.Schwartz}\
        The Schwartz space of $\GR$ with respect to $\phi$, $\SS(\GR,
        \phi)$, is closed under holomorphic functional calculus in
        $\cred{\GR}$.
\end{theorem}

\begin{proof}\
The methods of \cite{lafforgue-these} extend without difficulty to our
case.

As above, let $k_0$ be such that $\int_{\GR_x}
\frac{1}{(1+\phi(g))^{k_0}}$ converges for any $x \in \g0$.

\begin{lemma}\label{lemme2}\ If $l \geq k \geq k_0$ then
$\SS^{l,d}(\GR, \phi)$ is stable under holomorphic calculus in
$\SS^{k,d}(\GR, \phi)$.  
\end{lemma} 

\begin{proof}\ 
To prove this, 
by Lemma \ref{lemlaf}, we will show that
$$
        \lim_{n\to \infty} \n{f^n}{k,d}^{\frac{1}{n}}=\lim_{n\to
        \infty} \n{f^n}{l,d}^{\frac{1}{n}}.
$$ 
Indeed, this implies that $\SS^{l,d}(\GR, \phi)$, which is a dense
subalgebra of $\SS^{k,d}(\GR, \phi)$, is also full.

Now, as above, we have
\begin{eqnarray*}
 &  & |v_1 \ldots v_i \cdot
f^n\cdot v_{i+1}  \ldots v_l(g)|\\
 &=& |\int_{g_1\ldots g_n=g} (v_1 \ldots v_i
\cdot f)(g_1) f(g_2)\ldots  f(g_{n-1}) (f\cdot v_{i+1}  \ldots v_l)(g_n)| \\
 &\leq& \sum_{i=1}^{n} \int_{\begin{subarray}{l}
g_1\ldots g_n=g\\
\phi(g_i)\geq  \frac{\phi(g)}{n}
\end{subarray}}
  |v_1 \ldots v_i \cdot f(g_1)| |f(g_2)|\ldots  |f(g_{n-1})| |f\cdot v_{i+1}  
\ldots v_l(g_n)|\\
 & \leq& \sum_{i=1}^{n} \n{f}{k,d}^{n-1} \n{f}{l,d}
\int_{\begin{subarray}{l}
g_1\ldots g_n=g\\
\phi(g_i)\geq \frac{\phi(g)}{n}
\end{subarray}}
 \frac{1}{(1+\phi(g_i))^l} \prod_{j\neq i}\frac{1}{(1+\phi(g_j))^k}\\
 & \leq&  \n{f}{k,d}^{n-1} \n{f}{l,d}
\frac{1}{(1+\frac{\phi(g)}{n})^{l}} \sum_{i=1}^{n} \int_{\begin{subarray}{l}
g_1\ldots g_n=g\\
\phi(g_i)\geq \frac{\phi(g)}{n}
\end{subarray}}
\prod_{j \neq i}\frac{1}{(1+\phi(g_j))^k} \\
\end{eqnarray*}

As the latter integrals are lower than $C$,
$$ |f^n(g)|(1+\phi(g))^l \leq  \n{f}{k,d}^{n-1} \n{f}{l,d}
n^{l}\cdot nC$$
which gives
$$\n{f^n}{l,d}^\frac{1}{n} \leq n^{\frac{1+l}{n}} C^{\frac{1}{n}} 
\n{f}{k,d}^{1-\frac{1}{n}} \n{f}{l,d}^{\frac{1}{n}}$$
thus
$$\lim_{n\to \infty} \n{f^n}{l,d}^\frac{1}{n} \leq \n{f}{k,d}.$$
Let's now apply this inequality to $f^m$, we get
$$\lim_{n\to \infty} \n{f^{mn}}{l,d}^\frac{1}{n} \leq \n{f^m}{k,d},$$
so that if $m \to \infty$, 
$$ \lim_{n\to \infty} \n{f^{n}}{l,d}^\frac{1}{n} =\lim_{m,n\to \infty} 
\n{f^{mn}}{l,d}^\frac{1}{mn} \leq 
\lim_{m\to \infty} \n{f^{m}}{k,d}^\frac{1}{m}. $$
Because of  $\n{f^{m}}{k,d}\leq \n{f^{m}}{l,d}$  we get
$\lim_{n\to \infty} \n{f^{n}}{l,d}^\frac{1}{n} =\lim_{n\to \infty} 
\n{f^{n}}{k,d}^\frac{1}{n}$.
    \end{proof}

\begin{lemma}\label{lemme3}\
      If $k \geq k_0$, one has
$$
        \SS^{k,d}(\GR, \phi) * C^*(\GR) * \SS^{k,d}(\GR, \phi) 
        \subset \SS^{0,d}(\GR, \phi).
$$
\end{lemma}

\begin{proof}\
If  $f_1, f_2 \in \SS^{k,d}(\GR, \phi)$ and $f \in C^*(\GR)$, then
$$ 
        |v_1 \ldots v_i \cdot f_1 * f * f_2 \cdot v_{i+1}  \ldots v_l(g)|
        = | \int_{\GR_{s(g)}} (v_1 \ldots v_i \cdot f_1) (g{g'}^{-1}) (f*f_2 \cdot
        v_{i+1}  \ldots v_l)(g')|. 
$$
But if we denote by $f_3$ the function $g'\mapsto (v_1 \ldots v_i \cdot f_1)
(g{g'}^{-1})$, then the Cauchy-Schwarz inequality gives
\begin{equation*}
        \begin{split}
          |v_1 \ldots v_i \cdot f_1 * f * f_2 \cdot v_{i+1}  \ldots v_l(g)|
        \leq & \n{f_3}{L^2(\GR_{s(g)})} \n{f*f_2 \cdot
        v_{i+1}  \ldots v_l}{L^2(\GR_{s(g)})}\\
        \leq & C\n{f_3}{k,d} \n{f}{\cred{\GR}} \n{f_2}{L^2(G_{s(g)})}\\
        \leq & C\n{f_3}{k,d} \n{f}{\cred{\GR}} C \n{f_2}{k,d},
        \end{split}
\end{equation*}
and hence $f_1 * f * f_2 \in \SS^{0,d}(\GR, \phi)$. 
\end{proof}

\begin{lemma}\
If $k\geq k_0$, then $\SS^{k,d}(\GR, \phi)$ is stable under holomorphic 
functional calculus in $\cred{\GR}$.
\end{lemma}

\begin{proof}\
If $f \in \SS^{k,d}(\GR, \phi)$, then the lemma \ref{lemme3} implies that
$$
        \n{f^n}{0,d} \leq C^2\n{f}{k,d} \n{f^{n-2}}{\cred{\GR}} \n{f}{k,d}
$$ 
thus 
$$ 
        \lim_{n\to \infty} \n{f^n}{0,d}^{1/n} =\lim_{n\to \infty}
        \n{f^n}{\cred{\GR}}^{1/n}
$$
(the inverse inequality is given by Proposition \ref{subalgebra}).
To prove this lemma,  it remains to show that 
$$
        \lim_{n\to \infty} \n{f^n}{0,d}^{1/n}=\lim_{n\to \infty}
        \n{f^n}{k,d}^{1/n}.
$$ 

If $f_1, f_2 \in \SS^{2k,d}(\GR, \phi)$, then
\begin{equation*}
\begin{split}
        |v_1 \ldots v_i \cdot f_1 * f_2 \cdot v_{i+1} \ldots v_l(g)|
        &\leq \int_{\GR_x} |v_1 \ldots v_i \cdot f_1(g{g' }^{-1})|
        |f_2 \cdot v_{i+1} \ldots v_l(g')|\\ &\leq
        \int_{\begin{subarray}{l} \GR_{s(g)}\\ \phi(g{g'}^{-1})\geq
        \frac{\phi(g)}{2} \end{subarray}}
        \frac{\n{f_1}{2k,d}\n{f_2}{0,d}}{(1+\phi(g{g'}^{-1}))^k(1+\phi(g))^k}\\
        &+ \int_{\begin{subarray}{l} \GR_{s(g)}\\ \phi(g')\geq
        \frac{\phi(g)}{2} \end{subarray}}
        \frac{\n{f_1}{0,d}\n{f_2}{2k,d}}{(1+\phi(g'))^k(1+\phi(g))^k}\\
        & \leq \frac{C}{(1+\phi(g))^k}( \n{f_1}{0,d}\n{f_2}{2k,d}+
        \n{f_1}{2k,d}\n{f_2}{0,d}) \end{split}
\end{equation*}
hence $\n{f_1f_2}{k,d} \leq C( \n{f_1}{0,d}\n{f_2}{2k,d}+ 
\n{f_1}{2k,d}\n{f_2}{0,d})$.

Applying this to $f_1=f_2=f$, we get
$$
        \n{f^{2n}}{k,d} \leq 2C\n{f}{0,d}\n{f}{2k,d}.
$$
But by Lemma \ref{lemme2}, 
$$
        \lim_{n\to \infty} \n{f^n}{2k,d}^{1/n}=\lim_{n\to \infty}
        \n{f^n}{k,d}^{1/n},
$$
so that 
$$
        \lim_{n\to \infty} \n{f^n}{k,d}^{1/n} \leq \lim_{n\to \infty}
        \n{f^n}{0,d}^{1/n},
$$
which gives an equality as the opposite inequality comes from the Remark 
\ref{remark}.
\end{proof}

We are ready now to complete the proof of Theorem \ref{theorem.Schwartz}.
The Schwartz space of $\GR$ with respect to $\phi$ is thus an intersection of
subalgebras of $\cred{\GR}$ which are stable under holomorphic functional
calculus, hence $\SS(\GR,\phi)$ is  stable under holomorphic functional
calculus in $\cred{\GR}$.
\end{proof}

\begin{corollary}\
    Let $\GR$ be a \cfg, and $\phi$ be a length function with polynomial
    growth. Denote by $\Psi^0_s(\GR)$ the sum of $\Psi^0(\GR)$ and of the
    Schwartz space of $\GR$ with respect to $\phi$, $\SS(\GR,\phi)$.
    Then $\Psi^0_s(\GR)$ is stable under holomorphic functional calculus.

Moreover, if $P \in \Psi^0_s(\GR)$ is Fredholm, then it admits a parametrix $Q
\in \Psi^0_s(\GR)$.
\end{corollary}

\begin{proof}\

Let us first show that $\Psi_s^0(\GR)$ is an algebra, which amounts to prove that if $k_1, k_2 \in I_c^{m,0}(\GR, \GR^{(0)})$ and $f \in \SS(\GR, \phi)$ then $k_1*f*k_2 \in \SS(\GR, \phi)$.

The Lie algebroid of $\GR$ being a $\cinfo$ fiber bundle, it is possible to choose a set $v_1, \ldots, v_N \in \mathcal{C}(A(\GR))$ such that for any $x \in \GR^{(0)}$, $(v_1(x), \ldots, v_n(x))$ generates $T_x \GR_x$. Now let 
$\Delta=\sum_{i=1}^{N} v_i^2$. This differential operator is of degree 2 and admits a parametrix. 

If $k_1, k_2 \in I_c^{m,0}(\GR, \GR^{(0)})$ then for any $l$ there exist $k'_1, k'_2 \in I_c^{m-2l,0}(\GR, \GR^{(0)})$ and $r_1, r_2 \in \ccinfo(\GR)$ such that $k_1=k'_1 * k_\Delta^l +r_1, k_2=k_\Delta^l *k'_2 +r_2$. But when $k \in I_c^{-N-1,0}(\GR, \GR^{(0)})$, it is given by an absolutely convergent oscillatory integral, thus it belongs to $\Cc(\GR)$. Hence if $l$ is large enough, $k'_1$ and $k'_2$ belong to $\SS(\GR)$. But
$$k_1 *f* k_2=k'_1 * k_\Delta^l *f * k_\Delta^l * k'_2 +k'_1 * k_\Delta^l *f *r_2 +r_1*f*k_\Delta^l *k'_2  +r_1 *f*r_2$$
so that $k_1 *f* k_2 \in \SS(\GR)$ since $k_\Delta^l *f*k_\Delta^l \in \SS(\GR)$ by definition and $\SS(\GR)$ is an algebra.

The first part of this result is then a direct consequence of Theorem \ref{propred}.

Consider the exact sequence
$$
        0 \to \ck \to \alg{\GR} \to \alg{\GR}/\ck \to 0.
$$

If $P \in \Psi^0_s(\GR)$ is Fredholm, then its image in $\alg{\GR}/\ck$ is
invertible, thus, Theorem~\ref{thm4.02} implies that its inverse lives in
$\Psi^0_s(\GR)/(\Psi^0_s(\GR)\cap \ck)$. This shows that there exists a
parametrix in $\Psi^0_s(\GR)$.
\end{proof}

As P. Piazza has kindly informed us, Schwartz-type spaces seem to be important also
for the questions related to higher-signatures that were considered in 
\cite{piazza1,piazza2}.

\subsection{Schwartz spaces associated to the $c_n$ groupoids}
The \cfg s defined in Section \ref{Sec.4} can be endowed with length
functions in order to define, for each of them, a Schwartz space, thus
to obtain an algebra of \pd\ operators closed under holomorphic
functional calculus. To keep notations simple, we only consider here
the case of a compact manifold with boundary.  In the case of the
\bcalc, we have
$$
        \Gamma(X)=\{ (x,y,\lambda) \in X \times X \times \rep{},
        \rho(x)=\lambda \rho(y) \}.
$$
Define $\phi(x,y,\lambda) = |\log(\lambda)|$. It was shown in \cite{mont98}
that this defines a length function.

Before considering the case of the  $c_n$-calculi, we need an easy  lemma:
\begin{lemma}
  Let $\GR$ and $\GR'$ be two homeomorphic \cfg s. Then any length function 
with polynomial growth on $\GR$ induces  such a function on $\GR'$.
\end{lemma}
\begin{proof}
  Denote by $f : \GR' \to \GR$ the homeomorphism, and assume that 
$\GR$ has a length function with polynomial growth, $\phi$. 
Let $\phi'=\phi\circ 
f : \GR' \to \R_+$.  This function is clearly a length function with polynomial 
growth, since $f$ is a homeomorphism, and the Haar system on $\GR'$ is induced 
by that on $\GR$.
\end{proof}

Recall from Section 4 that there is a homeomorphism $\Theta_n:\Gamma(X) \to
\Gamma_{n}(X)$ where $\Gamma(X)$ is the \gr\ of the $b$-calculus, and
$\Gamma_{n}(X)$ is the \gr\ of the $c_n$-calculus.
\begin{proposition}
  Let $\phi_n=\phi \circ \Theta_n^{-1}$. Then $\phi_n$ is a length
  function. Under the identification 
$$
        \Gamma_{n+1}(X)=\{ (u,v,\mu) \in X\times X\times \R, 
        \mu \rho(u)^n\rho(v)^n=\rho(u)^n-\rho(v)^n\},
$$
one gets $\phi_n(u,v,\mu)=|\mu|$.
\end{proposition}
\begin{proof}
  This is clear since $\Theta_n(x,y,\lambda)=(u,v,\log(\lambda))$. 
\end{proof}


\end{document}